\documentclass{article}

\usepackage[german, english]{babel}
\usepackage{amsmath}
\usepackage{tabularray}
\usepackage{amsthm}
\usepackage{amssymb}
\usepackage{amsfonts,epsfig}
\usepackage{hyperref}
\usepackage{varwidth}
\usepackage{pifont}
\usepackage{blindtext}
\usepackage{upgreek}  
\usepackage[titletoc,title]{appendix}
\usepackage[bottom]{footmisc}
\usepackage[affil-it]{authblk}
\usepackage{enumerate}
\usepackage[shortlabels]{enumitem}
\usepackage{tikz}
\usetikzlibrary{shapes.geometric, arrows.meta,positioning,shapes.symbols,shapes.misc}
\usepackage{tikz-3dplot}
\usepackage{float}
\usepackage{version}
\usepackage{extarrows}
\usepackage{enumitem}
\usepackage{appendix}
\usepackage{mathrsfs}
\usepackage[scr=rsfs,cal=boondox]{mathalfa}
\usepackage{gauss}
\usepackage{graphicx}
\usepackage{caption,subcaption}
\usepackage{etoolbox}
\usepackage[margin=1.1in]{geometry}
\usepackage[T1]{fontenc}
\usepackage{array} 
\usepackage{nicematrix}
\NiceMatrixOptions{
	cell-space-top-limit = 4pt,
	cell-space-bottom-limit = 4pt
}
\usepackage[utf8]{inputenc}
\usepackage{amsmath, amsfonts, amssymb, algorithm, algpseudocode}
\setlist[enumerate]{itemsep=0ex}
\usepackage{chngcntr}
\usepackage{natbib}
\usepackage{calligra}

\usepackage{booktabs}

\renewcommand{\Im}{\operatorname{Im}}
\let\Im\undefined
\DeclareMathOperator{\Im}{Im}

\newcommand{\bmat}[1]{\begin{bmatrix}#1\end{bmatrix}}

\DeclareFontFamily{T1}{calligra}{}
\DeclareFontShape{T1}{calligra}{m}{n}{<->s*[1]callig15}{}
\DeclareMathAlphabet\mathcalligra   {T1}{calligra} {m} {n}
\DeclareMathAlphabet\mathzapf       {T1}{pzc} {mb} {it}
\DeclareMathAlphabet\mathchorus     {T1}{qzc} {m} {n}
\DeclareMathAlphabet\mathrsfso      {U}{rsfso}{m}{n}

\newtheoremstyle{theorem}
{10pt}		  
{10pt}  
{\sl}  
{\parindent}     
{\bf}  
{. }    
{ }    
{}     
\theoremstyle{theorem}
\newtheorem{theorem}{Theorem}
\newtheorem{corollary}{Corollary}
\newtheorem{proposition}{Proposition}
\newtheorem{remark}{Remark}
\newtheorem{example}{Example}
\newtheorem{lemma}{Lemma}
\newtheorem{definition}{Definition}
\newtheorem*{algorithmbox}{Algorithm}
\newtheorem{property}{Property}

\newcommand{\mycalli}[1]{{\calligra #1}}

\newcommand{\naturals}{\mathbb{N}}

\newcommand{\Field}{\mathbb{F}}

\newcommand{\reals}{\mathbb{R}}

\newcommand{\Res}{\mathfrak{R}}
\newcommand{\calr}{\mathcal{r}}

\newcommand{\nul}{\boldsymbol{\rm nul}}
\newcommand{\lnul}{\boldsymbol{\rm nul_\ell}}
\newcommand{\rank}{\boldsymbol{\rm rank}}
\newcommand{\codim}{\boldsymbol{\rm codim}}
\newcommand{\rrank}{\boldsymbol{\rm rank_r}}
\newcommand{\crank}{\boldsymbol{\rm rank_c}}
\renewcommand{\dim}{\boldsymbol{\rm dim}}
\newcommand{\card}{\boldsymbol{\rm card}}

\newcommand{\RS}{{\rm RS}}
\newcommand{\LNS}{{\rm LNS}}
\newcommand{\NS}{{\rm NS}}
\newcommand{\Span}{{\rm span}}

\newcommand{\CS}{{\rm CS}}

\newcommand{\URREF}{{\rm \mathbf{URREF}}}
\newcommand{\UREF}{{\rm \mathbf{UREF}}}
\newcommand{\LREF}{{\rm \mathbf{LREF}}}
\newcommand{\LRREF}{{\rm \mathbf{LRREF}}}

\newcommand{\Ub}{{\rm \mathbf{U}}}
\newcommand{\V}{{\rm \mathbf{V}}}
\newcommand{\Lb}{{\rm \mathbf{L}}}

\newcommand{\N}{{\rm \mathbf{N}}}
\newcommand{\M}{{\rm \mathbf{M}}}
\newcommand{\A}{{\rm \mathbf{A}}}
\newcommand{\B}{{\rm \mathbf{B}}}

\newcommand{\C}{{\rm \mathbf{C}}}

\newcommand{\e}{{\rm \mathbf{e}}}
\newcommand{\I}{{\rm \mathbf{I}}}
\newcommand{\R}{{\rm \mathbf{R}}}
\newcommand{\x}{{\rm \mathbf{x}}}
\newcommand{\Chi}{{\cal{X}}}

\newcommand{\y}{{\rm \mathbf{y}}}
\newcommand{\z}{{\rm \mathbf{z}}}
\newcommand{\n}{{\rm \mathbf{n}}}

\newcommand{\chib}{\boldsymbol{\chi}}
\newcommand{\Lambdab}{\boldsymbol{\mathsf{L}}}
\newcommand{\Usb}{\boldsymbol{\mathsf{U}}}

\newcommand{\res}{{\rm \mathbf{r}}}

\newcommand{\0}{{\rm \mathbf{0}}}

\newcommand{\Exnoether}{\mathcalligra{N\ }}

\newcommand{\Per}{{\boldsymbol{\rm P}}}
\newcommand{\Qb}{{\boldsymbol{\rm Q}}}
\newcommand{\Vb}{{\boldsymbol{\rm V}}}

\newcommand{\REQ}{\mathrel{\underset{\sim}{R}}\vspace{0.02in}}
\newcommand{\LEQ}{\mathrel{\underset{\sim}{L}}\vspace{0.02in}}

\makeatletter
\patchcmd\g@matrix
{\vbox\bgroup}
{\vbox\bgroup\normalbaselines}
{}{}
\makeatother

\newcommand{\BAR}[1][0pt]{%
	\hspace{-\arraycolsep}%
	\ifdim#1>0pt
	\dimen0=\ht\strutbox \dimen2=\dimen0
	\advance\dimen0 #1\relax
	\ht\strutbox=\dimen0
	\fi
	\smash{\strut\vrule} 
	\ifdim#1>0pt
	\ht\strutbox=\dimen2
	\fi
	\hspace{-\arraycolsep}%
}

\newcommand\utimes{\mathbin{\ooalign{$\cup$\cr%
			\hfil\raise0.42ex\hbox{$\scriptscriptstyle\times$}\hfil\cr}}}
\newcommand\bigutimes{\mathop{\ooalign{$\bigcup$\cr%
			\hfil\raise0.36ex\hbox{$\scriptscriptstyle\boldsymbol{\times}$}\hfil\cr}}}
\makeatletter		
\newcommand*{\pmzerodot}{%
	\nfss@text{%
		\sbox0{$\vcenter{}$}
		\sbox2{0}%
		\sbox4{0\/}%
		\ooalign{%
			0\cr
			\hidewidth
			\kern\dimexpr\wd4-\wd2\relax 
			\raise\dimexpr(\ht2-\dp2)/2-\ht0\relax\hbox{%
				\if b\expandafter\@car\f@series\@nil\relax
				\mathversion{bold}%
				\fi
				$\cdot\m@th$%
			}%
			\hidewidth
			\cr
			\vphantom{0}
		}%
	}%
}	

\makeatother

\begin{document}
	
	\bibliographystyle{plain}
	
\title{\textbf{\large 
The Canonical Family of Reduced Row Pivot Echelon Forms and the Induced Lexicographic Canonical Algorithm for  Linear Programming}
}	
\author{A. G. Paraskevopoulos$^{\dagger}$\thanks{Corresponding author: alparas@unipi.gr; paraskevopoulosalexandros@gmail.com} and D. A. Patsouras$^{\triangleright\star}$ \\
$^{\star}$University of Piraeus \ \  $^{\triangleright}$Costeas-Geitonas School.
\thanks{$^{\triangleright}$ The first author conceived the original concept and developed the formal proofs and algebraic framework. The second author proposed the finite-dimensional investigation and contributed insights leading to the Simplex--ILCA methodology. Both authors contributed to the conceptual refinement, preparation, and final revision of the manuscript.
}}
		
		\maketitle
		
\begin{abstract} {\footnotesize\noindent A complete canonical family of uniquely normalized, row-equivalent representatives of an $m\times n$ matrix is introduced and termed the family of Reduced Row Pivot Echelon Forms (RRPEFs). All RRPEFs represent the same row space, with the standard reduced row echelon form, renamed the Upper RREF (URREF), and the Lower RREF (LRREF) as its extreme members. A reformulation of Noether's Complementary Subspace Basis Theorem yields normalized bases of the row space that coincide with the nonzero rows of the corresponding RRPEFs. Constructively, each RRPEF is obtained by multi-index Gauss--Jordan reduction of a nonsingular transition submatrix. Factoring the admissible pivot-index sets by equality of their associated RRPEFs yields a quotient set in canonical bijection with the set of distinct RRPEFs. This bijection transfers the strict lexicographic order to the distinct RRPEFs and induces a finite chain from the URREF to the LRREF. The resulting structure gives rise to the Induced Lexicographic Canonical Algorithm (ILCA). Lexicographic confinement and singular-transition skipping restrict the candidate domain, whereas RRPEF reuse and scalar-equivalent transition bypassing avoid redundant Gauss--Jordan reductions while retaining the corresponding basic-solution tableaux. Strict lexicographic progression prevents cycling, while the traversal proceeds independently of objective improvement under degeneracy-induced stalling. ILCA terminates when a feasible tableau satisfies the standard optimality conditions and otherwise completes the finite traversal, including for Tucker- and Cunningham-type totally degenerate configurations. } \end{abstract}

\bigskip	
		
\noindent\textbf{Keywords:} Reduced Row Pivot Echelon Forms; canonical forms; Noether bases; Gauss--Jordan elimination; row equivalence; linear programming; pivoting algorithms; optimization.

 \bigskip 

\noindent\textbf{MSC 2020:} 15A21 (primary); 15A06, 90C05, 90C49 (secondary).

\section{Introduction} \label{sec.Introduction} In standard linear algebra, the RREF is recognized as the unique row-canonical representative of a finite $m\times n$ matrix under row equivalence. This paper extends the echelon paradigm by introducing a unified framework of multiple row-normal forms within the same row-equivalence class. To emphasize the structural symmetry of this framework, the standard RREF is termed the \emph{Upper} RREF (URREF), while its dual counterpart, the \emph{Lower} RREF (LRREF), is introduced in Definition~\ref{def. LREF} as the normal representative associated with the opposite pivot orientation. The LRREF arranges the nonzero rows according to their rightmost pivot positions, whereas the URREF uses the corresponding leftmost pivot positions. The coexistence of the URREF and LRREF may appear to conflict with the conventional uniqueness of a canonical form. Definition~\ref{def. RRPEF} resolves this apparent ambiguity by treating the individual Reduced Row Pivot Echelon Forms (RRPEFs) as row-normal representatives and transferring canonicality from the individual representatives to their complete family, which is uniquely determined by the row space. The RRPEF framework thereby encompasses the URREF and LRREF, defined in Definitions~\ref{def.UREF} and~\ref{def. LREF}, respectively, within a common axiomatic structure. Distinct RRPEFs generally differ both structurally and numerically, and transitions between them may require all three elementary row operations rather than row or column permutations alone; see Examples~\ref{ex. Noether basis} and~\ref{ex. WorkedExampleForRRPEFs}. Nevertheless, they span the same row space and belong to the same row-equivalence class.

\par A central contribution of this work is a reformulation of Noether's Complementary Subspace Basis Theorem (NCSBT; see~\cite{Ja:Ab}), which underlies the construction of multiple normalized bases of the row space. Motivated by the LRREF, the theorem is expressed through partitions of the standard basis containing admissible subsets that span complementary subspaces to the row space. A determinant test identifies the selections of standard unit vectors that constitute complementary bases, yielding the \emph{complementary criterion}. Specifically, when a row-space basis, such as the nonzero rows of an echelon form, is augmented by a selected set of standard unit vectors to form a square matrix, the determinant is nonzero precisely when the selected vectors form a basis of a complementary subspace. Each unit vector in the remaining part of the standard basis then admits a unique direct-sum decomposition into components in the row space and its complement. The resulting row-space components constitute a Noether basis whose elements coincide with the nonzero rows of the associated RRPEF; see Theorem~\ref{theo. normalized RRPEF}. Each admissible selection specifies the pivot columns relative to which the corresponding RRPEF is constructed. 

\par This construction yields row-equivalent normal forms extending from the URREF to the LRREF, as illustrated in Example~\ref{ex. Noether basis}, where complementary subspace bases are verified through the determinant test. Definition~\ref{def. RRPEF} unifies these forms within a single axiomatic framework. Because distinct Noether bases may be selected, finite and row-finite matrices may admit multiple normal representatives within the same row-equivalence class. In contrast to the finite case, where upper and lower row-normal forms coexist, this duality breaks down for row-finite matrices: only the LRREF provides a consistent row-equivalent extension. In infinite dimensions, the URREF fails to preserve row equivalence, as shown in~\cite{Pa:igjea2} and further analyzed in~\cite[Appendix~B]{Pa:LRREFROWFINITE}. Rooted in Toeplitz's seminal extension of finite-dimensional linear algebra to infinite matrices~\cite{Toeplitz1909}, and further developed through Fulkerson's quasi-Hermite form (QHF) for $\omega\times\omega$ row-finite matrices~\cite{Fulkerson1951}, this framework provides a unified methodology for reducing infinite matrices. The works in~\cite{Pa:igjea, Pa:reap, Pa:igjea2} address linear equations on infinite-dimensional spaces and row-finite systems, enabling explicit representations of solutions to linear difference equations with variable coefficients~\cite{Pa:Green}, with applications in time-series modeling; see~\cite{Ka:Pa, KaMRV14}. 

\par The complementary-basis index sets correspond bijectively, through set complementation, to the admissible pivot-column selections, generating the combinatorial index set $\mathcal{S}_{\mathbf{A}}$. Each element of this set identifies a nonsingular transition submatrix whose reduction to $\I_{\rho}$, where $\rho=\rank(\A)$, produces the associated RRPEF. The admissible pivot-index sets inherit the strict lexicographic order of $\binom{[n]}{\rho}$. Factoring $\mathcal{S}_{\mathbf{A}}$ by equality of the associated RRPEFs yields the quotient set $\mathcal{S}^o_{\mathbf{A}}$, termed the \emph{core combinatorial index set}. This quotient is identified bijectively with the set of minimal representatives of its equivalence classes and, in turn, with the set of distinct RRPEFs ${\rm Im}(\mathcal R_{\A})$; see Section~\ref{subsec:EquivalenceRelationOnTheCombSET}. Through this bijection, the strict lexicographic order $\prec_L$ transfers from the minimal representatives identified with $\mathcal{S}^o_{\mathbf{A}}$ to ${\rm Im}(\mathcal R_{\A})$, thereby inducing a finite RRPEF chain from the URREF to the LRREF; see Section~\ref{sec. The Hierarchy of normalized Forms}. 
\par This ordered algebraic structure gives rise to the \emph{Induced Lexicographic Canonical Algorithm} (ILCA), which traverses the admissible pivot-index sets in strict lexicographic order and obtains their associated RRPEFs. The invariance of the row space fixes the number of nonzero rows at $\rho$ in every row-equivalent normal form and confines the combinatorial domain to the single layer $\binom{[n]}{\rho}$. ILCA contracts the required computation through four mechanisms: lexicographic confinement, singular-transition skipping, RRPEF reuse, and scalar-equivalent transition bypassing. The first two restrict the candidate domain, whereas the latter two avoid redundant Gauss--Jordan reductions. In the linear-programming framework, each admissible pivot-index set determines an indexed basic-solution tableau whose constraint matrix is the associated RRPEF. Distinct admissible pivot-index sets may generate the same RRPEF while yielding different basic solutions and transformed objective rows. Accordingly, RRPEF reuse avoids repeated construction of identical RRPEFs without suppressing the formation and examination of the corresponding basic-solution tableaux. This establishes the structural link between the row-space algebra and the LP tableau framework while distinguishing the exhaustive traversal of admissible BCI sets from the construction of distinct RRPEFs. Although linear inequalities define a convex polyhedral feasible region, the construction and traversal of the RRPEF family are entirely algebraic and do not require feasibility-preserving pivoting. 

\par The Beale, Tucker, and Cunningham instances examined in Section~\ref{sec:degenerate} illustrate complementary aspects of the ILCA traversal. Cunningham's totally degenerate stalling example shows how lexicographic confinement and singular-transition skipping restrict the candidate domain, while RRPEF reuse avoids redundant Gauss--Jordan reductions without suppressing the corresponding basic-solution tableaux. These mechanisms operate independently of objective improvement and local pivot decisions. Strict lexicographic progression precludes the repetition of an indexed tableau and hence excludes cycling, while the finiteness of the BCI chain ensures termination despite degeneracy-induced stalling. In linear-programming applications, ILCA terminates when a feasible tableau satisfies the optimality conditions and otherwise completes the finite traversal, including on totally degenerate faces; see Property~\ref{prop. halting_rule}.

\section{Preliminaries and Notation}\label{sec. Preliminaries} This section collects the notation and basic matrix concepts used throughout the paper. We denote the sets of nonnegative and positive integers by $\naturals=\{0,1,2,\ldots\}$ and $\naturals_1=\{1,2,3,\ldots\}$, respectively. For $n\in\naturals_1$, let $\I(n)\triangleq\{1,\ldots,n\}$. For an indexed family $A=(A_i)_{i\in I}$, regarded as a function on the index set $I$, its image is $\Im(A)=\{A_i:i\in I\}$, the set of distinct values attained by the family. All vector spaces are defined over a field $\Field$ of characteristic zero. Matrices and vectors are written in \textbf{boldface}, scalars in \emph{italic}, and bases in {\Large \mycalli{calligraphic}}. For subspaces $G\sqsubset Y$, a complementary subspace is denoted by $G^c$, whereas the complement of a set $X$ is denoted by $X'$ and the transpose of a matrix $\A$ by $\A^{\top}$. To represent bases as indexed families, we use ordered indexed sets $\Chi=\{\chib_i\in Y\mid i\in\I(n)\}$. If $\Chi$ is a basis of $Y$, then $\dim(Y)=\card(\Chi)$. Let $\mathrsfso{E}_m=\{\e_i\mid i\in\I(m)\}$ denote the standard basis of the coordinate space $\Field^m$. For any matrix $\mathbf{A}\in\Field^{m\times n}$, its entries are denoted by $a_{ij}$, its $k$--th row by $\mathbf{A}_{k\bullet}\in\Field^{1\times n}$, and its $k$--th column by $\mathbf{A}_{\bullet k}\in\Field^{m\times 1}$. Under these conventions, the four fundamental subspaces are: 
\begin{itemize} 
\item $\CS(\mathbf{A})=\Span\{\mathbf{A}_{\bullet j}\mid j\in\I(n)\}\sqsubset\Field^{m\times 1}$ (column space),\vspace{-0.05in} 
\item $\RS(\mathbf{A})=\Span\{\mathbf{A}_{i\bullet}\mid i\in\I(m)\}\sqsubset\Field^{1\times n}$ (row space),\vspace{-0.05in} \item $\NS(\mathbf{A})=\{\mathbf{y}\in\Field^{n\times 1}\mid\mathbf{A}\mathbf{y}=\mathbf{0}\}$ (right null space),\vspace{-0.05in} \item $\LNS(\mathbf{A})=\{\mathbf{x}\in\Field^{1\times m}\mid\mathbf{x}\mathbf{A}=\mathbf{0}\}$ (left null space). \end{itemize} 
Their dimensions are denoted by the column rank $\crank(\mathbf{A})$, row rank $\rrank(\mathbf{A})$, nullity $\nul(\mathbf{A})$, and left (or row) nullity $\lnul(\mathbf{A})$, respectively. In finite dimensions, since 
$\rank(\mathbf{A})=\rrank(\mathbf{A})=\crank(\mathbf{A}),$
 the rank--nullity relations give\vspace{-0.1in}  \begin{equation}\label{equations rank2} \rank(\mathbf{A})=m-\lnul(\mathbf{A}),\qquad \rank(\mathbf{A})=n-\nul(\mathbf{A}).\vspace{-0.05in} \end{equation} 
Matrices $\mathbf{A},\mathbf{B}\in\Field^{m\times n}$ are \emph{row equivalent}, denoted by $\mathbf{A}\REQ\mathbf{B}$, if $\mathbf{B}$ is obtained from $\mathbf{A}$ by finitely many elementary row operations. They are \emph{left equivalent} if $\mathbf{B}=\mathbf{Q}\mathbf{A}$ for some invertible $\mathbf{Q}\in\Field^{m\times m}$ termed \emph{transition matrix}. Left and row equivalence relations are equivalent, since every finite product of elementary matrices is invertible and every invertible matrix is a finite product of elementary matrices. \begin{proposition}\label{prop. equivalent statements to row equivalence} Let $\mathbf{A},\mathbf{B}\in\Field^{m\times n}$. The following statements are equivalent:\vspace{-0.05in} \begin{enumerate}[\itshape i)] \item $\mathbf{A}\REQ\mathbf{B}$.\vspace{-0.05in} \item $\mathbf{A}\LEQ\mathbf{B}$.\vspace{-0.05in} \item $\RS(\mathbf{A})=\RS(\mathbf{B})$.\vspace{-0.02in} \item $\URREF(\mathbf{A})=\URREF(\mathbf{B})$. \end{enumerate} \end{proposition} 
\begin{definition}\label{def.UREF} {\rm A matrix $\Usb\in\Field^{m\times n}$ is in upper row echelon form (UREF) if: \begin{enumerate}[\itshape i)] \item The leftmost nonzero entry of each nonzero row is $1$, called its \emph{leading one}. \item The leading $1$ of each nonzero row lies strictly to the right of the leading $1$s of the preceding rows. \item All zero rows occur below the nonzero rows. \end{enumerate} 
A matrix $\Ub$ is in upper reduced row echelon form (URREF) if it is in UREF and additionally satisfies: \begin{enumerate}[\itshape iv)] 
\item Each pivot column contains a unique nonzero entry, namely its leading $1$. \end{enumerate} 
The URREF is the classical RREF.} 
\end{definition} 
If $\mathbf{U}$ satisfies Definition~\ref{def.UREF} and $\mathbf{A}\REQ\mathbf{U}$, then $\mathbf{U}=\URREF(\mathbf{A})$, the unique RREF of $\mathbf{A}$. Throughout, we write $\rho\triangleq\rank(\mathbf{A})$. The set $\mathsf{U}=\{\mathbf{U}_{1\bullet},\ldots,\mathbf{U}_{\rho\bullet}\}$ of nonzero rows of $\URREF(\mathbf{A})$ forms a basis of $\RS(\mathbf{A})$; hence $\card(\mathsf{U})=\rho=\dim(\RS(\mathbf{A}))$. \begin{proposition}\label{prop. UREF rank and row nullity} The following statements hold: \begin{enumerate}[\itshape i)] \item The nonzero rows of $\URREF(\mathbf{A})$ form a basis of $\RS(\A)$, and their number equals $\rank(\mathbf{A})$, which also equals the dimension of $\CS(\A)$. \item The number of zero rows of $\URREF(\mathbf{A})$ equals the row nullity, $\lnul(\mathbf{A})=m-\rho$. \end{enumerate} \end{proposition} The pivot columns of $\mathbf{U}=\URREF(\mathbf{A})$ are $\rho$ unit vectors $\mathbf{e}_i\in\mathrsfso{E}_m$ and form a basis of $\CS(\mathbf{U})$. Let $\mathcal{c}_\A$ denote the cardinality of this column basis. Then 
\[\mathcal{c}_\A=\dim(\CS(\mathbf{A}))=\dim(\RS(\A))=\rank(\mathbf{A})=\rho.\vspace{-0.02in} \] 
We shall also use the classical characterization of rank by nonzero minors . 
\begin{theorem}[Rank Minor Theorem]\label{theo. Rank Minor Theorem} Let \(\mathbf{A}\in\Field^{m\times n}\) and \(\rho\in\mathbb{N}\). The following statements are equivalent:\vspace{-0.05in}  
	\begin{enumerate}[\itshape i)] \item \(\operatorname{rank}(\mathbf{A})=\rho\). \item \(\mathbf{A}\) contains a \(\rho\times\rho\) submatrix \(\mathbf{M}\) with \(\det(\mathbf{M})\neq0\), and every \(k\times k\) submatrix with \(k>\rho\) has zero determinant. \end{enumerate} 
	\end{theorem} 
	By Rank Minor Theorem (RMT), the rank of a matrix is the order of its largest nonzero minor. The lower counterparts of UREF and URREF are defined by reversing the pivot orientation. 
	\begin{definition}\label{def. LREF} $\Lambdab\in\Field^{m\times n}$ is in \emph{LREF} if: \begin{enumerate}[\itshape i)] 
	\item The rightmost nonzero entry of each nonzero row is $1$, called its \emph{leading (or rightmost) one}.\vspace{-0.05in} 
	\item The leading $1$ of each nonzero row lies strictly to the right of the leading $1$s of the preceding rows.\vspace{-0.05in} 
	\item All zero rows occur below the nonzero rows.\vspace{-0.05in}
	 \end{enumerate} A matrix $\Lb$ is in \emph{LRREF} if it is in \emph{LREF} and additionally satisfies:\vspace{-0.07in} \begin{enumerate}[\itshape iv)] 
	\item Each pivot column contains a unique nonzero entry, namely its leading $1$.\vspace{-0.01in} \end{enumerate} \end{definition} \noindent Eqs.~\eqref{eq. 5x6 matrix} and~\eqref{equation LRREF of 5x6 matrix} illustrate a matrix $\A$ in URREF and its row-equivalent matrix $\Lb$ in LRREF. 
The UREF and LREF differ both numerically and structurally: their leading ones are, respectively, the leftmost and rightmost nonzero entries of their nonzero rows. Accordingly, Gauss--Jordan reduction under the leftmost-pivot convention produces the URREF, whereas Gauss--Jordan reduction under the rightmost-pivot convention produces the LRREF. These opposite pivot orientations are encoded in statement~{\rm (i)} of Definitions~\ref{def.UREF} and~\ref{def. LREF}. The URREF and LRREF are therefore row-canonical forms uniquely determined by their respective pivot conventions and constitute the two extreme members of the RRPEF family. More generally, each admissible pivot-index set determines a unique row-normal representative, whereas the complete RRPEF family is the canonical object uniquely determined by the row space.

\section{Noether Bases and Coordinate-Isolation Structures} \label{sec.Noether} This section develops the linear-algebraic foundations underlying the construction of the Reduced Row Pivot Echelon Forms (RRPEFs), including their two extreme members, the URREF and LRREF. The framework is based on complementary subspaces and their associated Noether bases. 

\subsection{Complementary Subspaces} \label{sub.Complementary} Let $Y$ be a finite- or infinite-dimensional vector space, and let $G$ be a subspace of $Y$. The subspaces $G$ and $G^{c}$ are \emph{complementary} if 
\[ G\cap G^{c} = \{\mathbf{0}\}, \qquad G+G^{c} = Y. \] Equivalently, 
\[ Y = G\oplus G^{c}; \] 
that is, every $\y\in Y$ admits a unique decomposition 
\[ \y = \y_G+\y_{G^{c}}, \qquad \y_G\in G, \qquad \y_{G^{c}}\in G^{c}. \] 
If $Y$ is infinite-dimensional, the existence of a complementary subspace for every subspace $G\sqsubset Y$ is guaranteed by Zorn's lemma. 
\par Suppose now that \(G\sqsubset Y\), with \(\dim(Y)=n<\infty\) and \(\dim(G)=m\le n\), and let \[ \mathcal{B}_{Y}=\{\mathbf{b}_1,\ldots,\mathbf{b}_n\}, \qquad \mathcal{B}_{G}=\{\mathbf{g}_1,\ldots,\mathbf{g}_m\} \] be bases of \(Y\) and \(G\), respectively. Then \[ \mathcal{B}_C=\{\mathbf{b}_{i_1},\ldots,\mathbf{b}_{i_{n-m}}\} \subseteq\mathcal{B}_Y \] is a basis of a complement \(G^c\) of \(G\) if and only if \(\mathcal{B}_G\cup\mathcal{B}_C\) is a basis of \(Y\). Different admissible selections $\mathcal{B}_C\subseteq\mathcal{B}_Y$ may yield different complementary subspaces of $G$, each of dimension \[ \dim(G^c) = \codim(G). \] 

\subsection{Noether's Complementary Basis Theorem}\label{sub.Noether} Let \(Y\) be a vector space, finite- or infinite-dimensional, and let \(G\) be a subspace of \(Y\). In Jacobson's formulation of the Noether Complementary Subspace Basis Theorem (NCSBT; \cite{Ja:Ab}, pp.~242--243), a basis \(\mathcal{B}_{Y}=\{\mathbf{b}_{i}\mid i\in I\}\) of \(Y\) is partitioned as 
\[ C=\{\mathbf{b}_{i}\mid i\in\mathsf{C}\}, \qquad B=\{\mathbf{b}_{i}\mid i\in\mathsf{C}'\}, \qquad \mathcal{B}_{Y}=C\cup B,\quad C\cap B=\emptyset, \] 
such that \(C\) is a basis of a complementary subspace \(G^{c}\) of \(G\). The corresponding \emph{Noether basis} of \(G\) is 
\[ \Exnoether_{\mathsf{C}} = \{\mathbf{n}_{i} \mid \mathbf{n}_{i} = \mathbf{b}_{i} + \mathbf{u}_{i}, \; \mathbf{u}_{i}\in G^{c}, \; i\in\mathsf{C}'\}. \] 
The NCSBT relies on the axiom of choice, invoked through Zorn's lemma, to guarantee the existence of a complementary subspace for a given subspace.~\footnote{An alternative version of NCSBT was established in \cite{Pa:reap} and adapted to row-finite matrices in \cite{Pa:LRREFROWFINITE} through well-ordered bases. This approach replaces Zorn's lemma with Zermelo's well-ordering theorem while remaining within ZFC, since both are equivalent to the axiom of choice. Thus, the scope and substance of NCSBT are preserved while facilitating the extension from finite row-canonical matrices to row-finite ones.} 
\begin{definition}[CBI Set]\label{def.CBI} Let \(\mathcal{B}_{Y}=\{\mathbf{b}_{i}\mid i\in I\}\) be a basis of \(Y\) and let \(\mathsf{C}\subset I\) be such that \[ C=\{\mathbf{b}_{i}\mid i\in\mathsf{C}\} \] is a basis of a complementary subspace \(G^{c}(\mathcal{B}_{Y},\mathsf{C})\) of \(G\). Then \(\mathsf{C}\) is called a \emph{complementary basis index set} (\emph{CBI} set) associated with \(G\). The collection of all CBI sets associated with \(G\) is denoted by \(\mathcal{C}_{G}\). 
\end{definition} 
For fixed \(\mathcal{B}_{Y}\) and \(\mathsf{C}\in\mathcal{C}_{G}\), the direct sum \[ Y=G\oplus G^{c}(\mathcal{B}_{Y},\mathsf{C}) \] yields, for each \(i\in\mathsf{C}'\), a unique decomposition, written for consistency with the Noether-basis representation as \[ \mathbf{b}_{i} = \mathbf{n}_{i}-\mathbf{u}_{i}, \qquad \mathbf{n}_{i}\in G,\quad \mathbf{u}_{i}\in G^{c}(\mathcal{B}_{Y},\mathsf{C}), \] 
or equivalently, 
\[ \mathbf{n}_{i}=\mathbf{b}_{i}+\mathbf{u}_{i}. \] Hence, each \(\mathbf{n}_{i}\), and therefore the corresponding Noether basis, is uniquely determined for each given pair $(\mathcal{B}_{Y},\mathsf{C})$. This yields the following reformulation of the NCSBT in terms of CBI sets.
\begin{theorem}[RNCSBT]\label{theo.RNCSBT} Let \(\mathcal{B}_{Y}=\{\mathbf{b}_{i}\mid i\in I\}\) be a basis of \(Y\), and let \(G\) be a subspace of \(Y\). For any CBI set \(\mathsf{C}\) associated with \(G\), the set \[ \Exnoether_{\mathsf{C}} = \{\mathbf{n}_{i}\mid i\in\mathsf{C}'\}, \qquad \mathbf{n}_{i} = \mathbf{b}_{i} + \mathbf{u}_{i}, \quad \mathbf{u}_{i}\in G^{c}(\mathcal{B}_{Y},\mathsf{C}), \] 
is the uniquely determined Noether basis of \(G\) relative to \((\mathcal{B}_{Y},\mathsf{C})\). 
\end{theorem} In finite dimensions, the RNCSBT yields a constructive procedure for determining all CBI sets through determinant-based admissibility tests. This construction is elementary and requires no appeal to the axiom of choice, whether through Zorn's lemma or Zermelo's well-ordering theorem.

\subsection{Complements of Row Spaces}\label{sub.ComplementaryRowSpace} Let 
\[ \A\in\Field^{m\times n}, \qquad \rank(\A)=\rho, \qquad \boldsymbol{\mathsf U}=\UREF(\A), \] 
and let 
\[ \mathsf U= \left\{ \boldsymbol{\mathsf U}_{1\bullet},\ldots, \boldsymbol{\mathsf U}_{\rho\bullet} \right\} \] 
e the set of nonzero rows of \(\boldsymbol{\mathsf U}\). Since \(\mathsf U\) is a basis of \(\RS(\A)\), 
\[ \card(\mathsf U)=\dim\big(\RS(\A)\big)=\rank(\A)=\rho. \] Applying the construction of Section~\ref{sub.Complementary} with 
\[ Y=\Field^n, \qquad G=\RS(\A), \qquad \mathcal B_G=\mathsf U, \qquad \mathcal B_Y=\mathrsfso E_n, \] reduces the determination of a complementary subspace of \(\RS(\A)\) to selecting \(n-\rho\) standard basis vectors that, together with \(\mathsf U\), form a basis of \(\Field^n\). \par Remove the \(m-\rho\) zero rows of \(\boldsymbol{\mathsf U}\) and append \(n-\rho\) selected standard basis vectors $\{\e_{i_1}\ldots \e_{i_{n-\rho}}\}$ to its \(\rho\) nonzero rows, forming 
\[ \widetilde{\boldsymbol{\mathsf U}}=\mathsf U\cup\{\e_{i_1}\ldots \e_{i_{n-\rho}}\}\in\Field^{n\times n}. \] 
The selected standard basis vectors form a basis of a complementary subspace of \(\RS(\A)\) if and only if \[ \det\bigl(\widetilde{\boldsymbol{\mathsf U}}\bigr)\neq0. \] Thus, each admissible selection determines a complementary subspace of \(\RS(\A)\). Its index set specifies the complementary coordinate positions and, by set complementation, the associated admissible pivot positions of the corresponding row-normal form.

 \subsection{Matrix Realization}\label{sub.MatrixImplementation} Applying Theorem~\ref{theo.RNCSBT} to \(G=\RS(\A)\), for \(\A\in\Field^{m\times n}\), define \[ \Res_{\A}(\mathrsfso{E}_{n},\mathsf{C}) \triangleq \RS(\A)^{c}(\mathrsfso{E}_{n},\mathsf{C}), \] called the \emph{residual subspace} of \(\A\) relative to \((\mathrsfso{E}_{n},\mathsf{C})\).~\footnote{The term ``residual space'' was introduced by Ortiz~\cite{Or:tau} in the context of linear differential operators with polynomial coefficients acting on spaces of polynomials. This complementary interpretation of residual spaces was generalized in~\cite{Pa:reap} to abstract linear operators.} A CBI set of \(\RS(\A)\) is called a CBI set associated with \(\A\), and the collection of all such sets is denoted by \(\mathcal{C}_{\A}\). For each \(\mathsf{C}\in\mathcal{C}_{\A}\), \[ \Field^{n} = \RS(\A)\oplus \Res_{\A}(\mathrsfso{E}_{n},\mathsf{C}). \] For each \(i\in\mathsf{C}'\), the corresponding Noether basis element is uniquely determined by \[ \mathbf{n}_{i} = \mathbf{e}_{i}+\res_{i}, \qquad \res_{i}\in \Res_{\A}(\mathrsfso{E}_{n},\mathsf{C}), \] and hence has the coordinate representation \begin{equation}\label{Noether basis 1} \n_i = \e_i+\sum_{q\in\mathsf{C}}r_{iq}\e_q, \qquad i\in\mathsf{C}'. \end{equation} This matrix realization of the RNCSBT is summarized as follows. \begin{proposition}\label{prop.FiniteMatrixNoether} For any CBI set \(\mathsf{C}\in\mathcal{C}_{\A}\), the indexed set \[ \Exnoether(\mathrsfso{E}_{n},\mathsf{C}) = \{\mathbf{n}_{i}\mid i\in\mathsf{C}'\}, \] with \(\mathbf{n}_{i}\) defined by~\eqref{Noether basis 1}, is the uniquely determined Noether basis of \(\RS(\A)\) relative to \((\mathrsfso{E}_{n},\mathsf{C})\). \end{proposition} In the general setting, the notation \(\Exnoether(\mathcal{B}_{Y},\mathsf{C})\) records the dependence of the Noether basis on both the ambient basis \(\mathcal{B}_{Y}\) and the CBI set \(\mathsf{C}\). In the matrix setting, the ambient basis is fixed as the standard basis \(\mathrsfso{E}_{n}\) of \(\Field^{1\times n}\), relative to which \[ \A_{i\bullet} = \e_i\A = \sum_{j=1}^{n}a_{ij}\e_j, \qquad i\in\I_m. \] Here, the unit vector \(\e_i\) in \(\e_i\A\) belongs to \(\mathrsfso{E}_{m}\), whereas the unit vectors \(\e_j\) in the coordinate expansion belong to \(\mathrsfso{E}_{n}\). Henceforth, the fixed ambient basis is suppressed from the notation: \[ \Exnoether(\mathrsfso{E}_{n},\mathsf{C}) \equiv \Exnoether_{\mathsf{C}}, \qquad \Res_{\A}(\mathrsfso{E}_{n},\mathsf{C}) \equiv \Res_{\A}(\mathsf{C}). \]
\begin{example}\label{ex. Noether basis} \rm Consider the subspace 
\[ G=\{[x,y,z]\in\mathbb{R}^{3}\mid x+y+z=0\}. \] The vectors \[ \mathbf{g}_{1}=[1,1,-2], \qquad \mathbf{g}_{2}=[1,-2,1] \] 
form a basis \(\mathcal{G}=\{\mathbf{g}_{1},\mathbf{g}_{2}\}\) of \(G\), and \(\codim(G)=1\). Let \(\mathrsfso{E}_{3}=\{\e_{1},\e_{2},\e_{3}\}\) be the standard basis of \(\mathbb{R}^{3}\). For each \(i\in\{1,2,3\}\), 
\[ D_i= \det \begin{bmatrix} \mathbf{g}_{1}\\ \mathbf{g}_{2}\\ \e_i \end{bmatrix} =-3\neq0. \] 
Hence, each \(\e_i\) spans a complementary subspace of \(G\), so that 
\[ \mathcal C_G=\bigl\{\{1\},\{2\},\{3\}\bigr\}. \] 
\par\medskip \noindent\textbf{(i) Case \(\mathsf{C}=\{3\}\).} Here 
\[ G^{c}(\mathrsfso{E}_{3},\{3\})=\Span\{\e_{3}\}, \] 
and the Noether basis elements have the form 
\[ \mathbf{n}_{1}=\e_{1}+\beta_{1}\e_{3}, \qquad \mathbf{n}_{2}=\e_{2}+\beta_{2}\e_{3}. \] 
Imposing \(\mathbf{n}_{1},\mathbf{n}_{2}\in G\) gives \(\beta_{1}=\beta_{2}=-1\), and hence 
\[ \mathbf{n}_{1}=[1,0,-1], \qquad \mathbf{n}_{2}=[0,1,-1]. \] 
Thus,
\[ \mathbf{U} = \begin{bmatrix} 1&0&-1\\ 0&1&-1 \end{bmatrix}. \] 
Let us call 
\[ \mathbf{G} = \begin{bmatrix} 1&1&-2\\ 1&-2&1 \end{bmatrix}. \] 
Since \[ \RS(\mathbf{G}) = G = \RS(\mathbf{U}), \] 
we have \(\mathbf{G}\REQ\mathbf{U}\); moreover, \(\mathbf U\) satisfies Definition~\ref{def.UREF}. Therefore, 
\[ \mathbf{U}=\URREF(\mathbf{G}). \] \par\medskip \noindent\textbf{(ii) Case \(\mathsf{C}=\{1\}\).} Since \(\mathsf{C}'=\{2,3\}\), the corresponding Noether basis is \[ \mathbf{n}_{2}=[-1,1,0], \qquad \mathbf{n}_{3}=[-1,0,1], \] yielding \[ \mathbf{L} = \begin{bmatrix} -1&1&0\\ -1&0&1 \end{bmatrix} = \LRREF(\mathbf{G}). \] \par\medskip \noindent\textbf{(iii) Case \(\mathsf{C}=\{2\}\).} Since \(\mathsf{C}'=\{1,3\}\), the corresponding Noether basis is \[ \mathbf{n}_{1}=[1,-1,0], \qquad \mathbf{n}_{3}=[0,-1,1], \] yielding \begin{equation}\label{eq. Matrix C in Exam. 2} \C= \begin{bmatrix} 1&-1&0\\ 0&-1&1 \end{bmatrix}. \end{equation} The matrix \(\C\) is neither in URREF nor in LRREF, since its prescribed pivot positions satisfy neither the leftmost- nor the rightmost-pivot convention of Definitions~\ref{def.UREF} and~\ref{def. LREF}. Nevertheless, it is the uniquely normalized row-equivalent representative associated with \(\mathsf{C}=\{2\}\), with complementary pivot-index set \(\mathsf{C}'=\{1,3\}\). This intermediate form motivates the RRPEF framework developed in the next section. \end{example}

\section{The RRPEF Framework}\label{sec. An Extended Framework for normalized Forms} In this section, we develop a unified framework for the complete family of normalized forms introduced in Definition~\ref{def. RRPEF}. The framework combines four principal components: (i) the identification of the Noether bases of the row space with the nonzero rows of the RRPEFs of a matrix; (ii) the complement duality between the CBI sets in $\mathcal{C}_{\A}$ and the admissible pivot-column index sets in $\mathcal{S}_{\A}$; (iii) the bijection between the quotient set $\mathcal{S}^o_{\A}$ and the image set of the indexed RRPEF family; and (iv) the completeness of this family. The CBI sets provide the complementary coordinate structure underlying the construction of the RRPEFs and their classification as minimal, non-extreme, or maximal forms. A column position is a pivot position relative to an RRPEF if it indexes a column containing a pivot $1$, whereas the complementary column positions are indexed by the corresponding CBI set.\footnote{In the coordinate-isolation framework for systems of linear equations, the pivot columns correspond to the basic variables indexed by the BCI sets $\mathsf{S}\in\mathcal{S}_{\A}$, whereas the complementary columns correspond to the free coordinates indexed by the CBI sets. The complement operator $\mathsf{C}\mapsto\mathsf{C}'=\mathsf{S}$ identifies these two index structures.} Standard Gauss--Jordan elimination preserves the pivot positions selected by Gaussian elimination. The resulting pivot columns correspond to the leftmost nonzero entries of the nonzero rows of $\UREF(\A)$ and are called the \emph{leftmost pivot positions} (LPPs) of $\A$; their index set is denoted by $\mathsf{S}_{\ell}$. Dually, the \emph{rightmost pivot positions} (RPPs) correspond to the rightmost nonzero entries of the nonzero rows of $\Lambdab\triangleq\LREF(\A)$; their index set is denoted by $\mathsf{S}_{\calr}$. The existence of the LRREF motivates a broader notion of normalization in which the pivot positions are not restricted to the leftmost nonzero entries. Instead, each admissible pivot-index set determines a unique normalized row-equivalent representative. \subsection{The BCI Collection $\mathcal{S}_{\A}$}\label{subsec: The Combinatorial Collection} In view of Theorem~\ref{theo. Rank Minor Theorem}, we define the combinatorial collection associated with $\A$ of rank $\rho$. The elements of this collection generalize the sets $\mathsf{S}_{\ell}$ and $\mathsf{S}_{\calr}$. Throughout the remainder of the paper, the collection of all $\rho$-element subsets of $\I(n)$ is denoted by \[ \binom{[n]}{\rho}\triangleq\left\{\mathsf{S}\subseteq\I(n)\mid\card(\mathsf{S})=\rho\right\}. \] Thus, $\binom{[n]}{\rho}$ consists of the $\binom{n}{\rho}$ subsets of cardinality $\rho$ selected from $\I(n)=\{1,\ldots,n\}$. 
\begin{definition}[Basis Columns Indexing Set]\label{def. basis_indexing_set} {\rm Let $\Ub\in\Field^{\rho\times n}$ be the full-row-rank nonzero block of $\UREF(\A)$. A set \[ \mathsf{S}=\{\mathcal{i}_1,\mathcal{i}_2,\ldots,\mathcal{i}_\rho\}\in\binom{[n]}{\rho},\qquad \mathcal{i}_1<\cdots<\mathcal{i}_\rho, \] is called a \emph{Basis Columns Indexing} (BCI) set of $\A$ if the transition submatrix \[ \Ub_{\mathsf{S}}\triangleq\Ub_{\mathsf I_\rho,\mathsf S}=\begin{bmatrix}\Ub_{\bullet\mathcal{i}_1}&\cdots&\Ub_{\bullet\mathcal{i}_\rho}\end{bmatrix}\in\Field^{\rho\times\rho} \] is nonsingular. The combinatorial collection of all BCI sets associated with $\A$ is defined by \[ \mathcal{S}_{\A}\triangleq\left\{\mathsf{S}\in\binom{[n]}{\rho}\ \middle|\ \det(\Ub_{\mathsf{S}})\neq0\right\}. \] } \end{definition}
\begin{remark}\label{rem. normalized and full rank} {\rm By the Rank Minor Theorem, every matrix $\A\in\Field^{m\times n}$ of rank $\rho$ contains a nonsingular $\rho\times\rho$ minor $\A_{\mathsf{J},\mathsf{S}}$, indexed by a row-index set $\mathsf{J}\subseteq\I(m)$ and a column-index set $\mathsf{S}\subseteq\I(n)$, both of cardinality $\rho$. In the normalized matrix $\Ub=\URREF(\A)$, the $\rho$ nonzero rows occupy the first $\rho$ positions. Hence, the row-index set is fixed as $\I(\rho)$, whereas the admissible column-index set $\mathsf{S}$ may vary. Since elementary row operations preserve the linear independence of the columns indexed by $\mathsf{S}$, the BCI collection is invariant under row equivalence. Accordingly, the $\rho\times\rho$ transition submatrix associated with $\mathsf{S}$ is denoted by 
	\[ \Ub_{\mathsf{S}}\triangleq\Ub_{\I(\rho),\mathsf{S}}, \] where the fixed row-index set $\I(\rho)$ is omitted from the notation. Thus, all combinatorial variability is confined to the BCI sets \[ \mathsf{S}\in\mathcal{S}_{\A}\subseteq\binom{[n]}{\rho}. \] } 
\end{remark}
 A set $\mathsf{S}\in\binom{[n]}{\rho}$ can therefore be tested for membership in $\mathcal{S}_{\A}$ through the normalized matrix $\Ub=\URREF(\A)$.\\
 
 \bigskip 
 
 \noindent\textbf{Column-Index Pivot Criterion}\\ For the transition submatrix $\Ub_{\mathsf{S}}\in\Field^{\rho\times\rho}$ defined above, 
\[ \mathsf{S}\in\mathcal{S}_{\A}\quad\Longleftrightarrow\quad\det(\Ub_{\mathsf{S}})\neq0. \]

 \subsection{Axiomatic Characterization of RRPEFs} \rm The following definition extends Definitions~\ref{def.UREF} and~\ref{def. LREF} by allowing any admissible pivot-index set in $\mathcal{S}_{\A}$, thereby unifying the normalized row-equivalent forms of $\A$. 
 \begin{definition}[Reduced Row Pivot-Echelon Form]\label{def. RRPEF} \normalfont Let $\A\in\Field^{m\times n}$ be of rank $\rho$. A matrix $\R_{\mathsf{S}}=[r_{jk}]\in\Field^{m\times n}$ is said to be a RRPEF of $\A$ relative to \[ \mathsf{S}=\{\mathcal{i}_1,\ldots,\mathcal{i}_{\rho}\}\in\mathcal{S}_{\A} \] 
 if $\R_{\mathsf{S}}$ is row equivalent to $\A$, that is, \[ \R_{\mathsf{S}}\REQ\A, \] and the following four postulates hold: \begin{enumerate}[\itshape i)] 
 \item For every $j\in\I(\rho)$, the $j$-th nonzero row $(\R_{\mathsf{S}})_{j\bullet}$ has a pivot equal to $1$ in position $(j,\mathcal{i}_j)$; that is, \[ r_{j\mathcal{i}_j}=1. \] \item Equip $\I(\rho)$ and $\mathsf{S}$ with the standard strict order $<$ on $\naturals$. The mapping 
 \begin{equation}\label{eq. order isomorphism} \mathcal{i}:\I(\rho)\longrightarrow\mathsf{S},\qquad j\longmapsto\mathcal{i}_j, \end{equation} 
 is an \emph{order isomorphism}; that is, for every $i,j\in\I(\rho)$, 
 \[ i<j\quad\Longleftrightarrow\quad\mathcal{i}_i<\mathcal{i}_j. \] 
 Equivalently, the pivot positions satisfy \[ \mathcal{i}_1<\cdots<\mathcal{i}_{\rho}, \] so the pivot in each nonzero row lies strictly to the right of the pivots in all preceding rows. \item Every zero row of $\R_{\mathsf{S}}$ is located at the bottom of the matrix. \item Every column containing a pivot $1$ has all its remaining entries equal to $0$. \end{enumerate}  \end{definition}
 
 \subsection{The Lexicographic Order on BCI Sets}\label{subsec: Lexicographic Order on BCI Sets} Let \[ \mathcal{A}=(a_1,\ldots,a_k)\qquad\text{and}\qquad\mathcal{D}=(d_1,\ldots,d_k) \] be strictly increasing finite sequences of equal length $k$. The strict lexicographic order $\prec_L$ is defined by \[ \mathcal{A}\prec_L\mathcal{D} \] if there exists an index $\mathcal{h}\in\I(k)$ such that \[ a_i=d_i\quad\text{for all }i<\mathcal{h},\qquad\text{and}\qquad a_{\mathcal{h}}<d_{\mathcal{h}}, \] where $<$ denotes the standard strict order on $\naturals$. Thus, the first position at which the two sequences differ determines their relative order. Moreover, \[ \mathcal{A}=\mathcal{D}\quad\Longleftrightarrow\quad\mathcal{A}\not\prec_L\mathcal{D}\ \text{and}\ \mathcal{D}\not\prec_L\mathcal{A}. \] Let $\A\in\Field^{m\times n}$ have rank $\rho$. The number of $\rho$-element candidate pivot-index sets is \begin{equation}\label{eq. cardinality candidate pivot sets} \binom{n}{\rho}=\frac{n!}{(n-\rho)!\,\rho!}, \end{equation} which is the cardinality of $\binom{[n]}{\rho}$. 
 Each set \[ \mathsf{S}=\{\mathcal{i}_1,\ldots,\mathcal{i}_{\rho}\}\in\binom{[n]}{\rho} \] has the unique strictly increasing array configuration
 \[(\mathcal{i}_1,\ldots,\mathcal{i}_{\rho}),\qquad\mathcal{i}_1<\cdots<\mathcal{i}_{\rho}. \] Throughout, however, we retain the set notation \[ \mathsf{S}=\{\mathcal{i}_1,\ldots,\mathcal{i}_{\rho}\}, \] with the increasing order of its elements understood.
 
 Accordingly, the strict lexicographic order on these representations induces a strict total order on $\binom{[n]}{\rho}$ and, by restriction, on $\mathcal{S}_{\A}$. The BCI and CBI sets may therefore be compared lexicographically through their unique strictly increasing representations, while retaining their set-theoretic terminology. 
 
 \subsection{The Quotient Set $\mathcal{S}_{\A}^{o}$}\label{subsec:EquivalenceRelationOnTheCombSET} In general, \[ \card(\mathcal{S}_{\A})\geq\card\bigl(\Im(\mathcal{R}_{\A})\bigr), \] because distinct BCI sets may generate the same RRPEF, as illustrated by Cunningham's example. For $\mathsf{S}_i,\mathsf{S}_j\in\mathcal{S}_{\A}$, define \[ \mathsf{S}_i\sim_{\mathcal{S}_{\A}}\mathsf{S}_j\quad\Longleftrightarrow\quad\R_{\mathsf{S}_i}=\R_{\mathsf{S}_j}. \] Equality of matrices is reflexive, symmetric, and transitive; hence $\sim_{\mathcal{S}_{\A}}$ is an equivalence relation on $\mathcal{S}_{\A}$. The corresponding quotient set is \[ \mathcal{S}_{\A}^{o}\triangleq\mathcal{S}_{\A}/\mathord{\sim_{\mathcal{S}_{\A}}}. \] Since $(\mathcal{S}_{\A},\prec_L)$ is strictly totally ordered, every equivalence class $[\mathsf{S}]\in\mathcal{S}_{\A}^{o}$, being a nonempty finite subset of $\mathcal{S}_{\A}$, has a unique least element under $\prec_L$. Define \[ \mathcal{M}\triangleq\left\{\min_{\prec_L}[\mathsf{S}]\mid[\mathsf{S}]\in\mathcal{S}_{\A}^{o}\right\}. \] The mapping \[ f:\mathcal{S}_{\A}^{o}\longrightarrow\mathcal{M},\qquad f([\mathsf{S}])=\min_{\prec_L}[\mathsf{S}], \] is bijective. Indeed, every equivalence class has a unique least element, and distinct equivalence classes have distinct least elements. Since $\mathcal{M}\subseteq\mathcal{S}_{\A}$, the strict lexicographic order $\prec_L$ restricts to a strict total order on $\mathcal{M}$. Through the bijection $f$, this order transfers uniquely to $\mathcal{S}_{\A}^{o}$, making $f$ an order isomorphism. We therefore identify the quotient set $\mathcal{S}_{\A}^{o}$ with its canonical representative realization $\mathcal{M}$. \begin{definition}[Core Combinatorial Index Set]\label{def. working core BCI set} {\rm Under the preceding identification, \[ \mathcal{S}_{\A}^{o}\equiv\mathcal{M}=\left\{\min_{\prec_L}[\mathsf{S}]\mid[\mathsf{S}]\in\mathcal{S}_{\A}/\mathord{\sim_{\mathcal{S}_{\A}}}\right\}. \] Equivalently, \[ \mathsf{S}\in\mathcal{S}_{\A}^{o}\quad\Longleftrightarrow\quad\mathsf{S}\in\mathcal{S}_{\A}\ \text{and}\ \mathsf{S}=\min_{\prec_L}[\mathsf{S}]. \] The set $\mathcal{S}_{\A}^{o}$ is called the \emph{core combinatorial index set} associated with $\A$. } \end{definition}
 
\subsection{Noether Bases and RRPEFs} The following proposition identifies the nonzero rows of an RRPEF with the Noether basis associated with its complementary CBI set. The converse construction is established in Theorem~\ref{theo. normalized RRPEF}, and together these results yield a bijection between the sets of distinct Noether bases and distinct RRPEFs. \begin{proposition}\label{prop. rows of RRPEF basis of RS} Let $\A\in\Field^{m\times n}$ have rank $\rho$, and let $\R_{\mathsf{S}}=[r_{jk}]\in\Field^{m\times n}$ be an \emph{RRPEF} of $\A$ relative to $\mathsf{S}\in\mathcal{S}_{\A}$. Then: \begin{enumerate}[\itshape i)] \item The complementary set \[ \mathsf{S}'\triangleq\I(n)\setminus\mathsf{S} \] is a \emph{CBI} set; that is, \[ \mathsf{S}'=\mathsf{C}\in\mathcal{C}_{\A}. \] \item The set \[ \mathcal{N}\triangleq\left\{(\R_{\mathsf{S}})_{j\bullet}\mid j\in\I(\rho)\right\} \] of nonzero rows of $\R_{\mathsf{S}}$ is the Noether basis of $\RS(\A)$ uniquely determined by $\mathsf{C}=\mathsf{S}'$. \end{enumerate} \end{proposition} \begin{proof} By Definition~\ref{def. RRPEF}, the columns of $\R_{\mathsf{S}}$ indexed by $\mathsf{S}$ contain a $\rho\times\rho$ identity submatrix. Hence, the $\rho$ nonzero rows of $\R_{\mathsf{S}}$ are linearly independent. Since $\R_{\mathsf{S}}\REQ\A$ and $\rank(\A)=\rho$, the set $\mathcal{N}$ forms a basis of $\RS(\A)$. 
	
\noindent\textit{i)} It suffices to show that the set of standard unit vectors \[ \mathcal{Q}\triangleq\{\e_q\mid q\in\mathsf{S}'\} \] spans a complementary subspace of $\RS(\A)$; that is, \[ \Span(\mathcal{Q})=\Res_{\A}(\mathsf{S}'). \] Append the vectors in $\mathcal{Q}$ below the nonzero rows in $\mathcal{N}$ to form the square $n\times n$ matrix \begin{equation}\label{eq. Appended Matrix} \widetilde{\R_{\mathsf{S}}}\triangleq\left[\begin{array}{c}\mathcal{N}\vspace{0.05in}\\ \hline\vspace{-0.1in}\\ \mathcal{Q}\end{array}\right]. \end{equation} The identity $\Span(\mathcal{Q})=\Res_{\A}(\mathsf{S}')$ follows once the rows of $\widetilde{\R_{\mathsf{S}}}$ are shown to form a basis of $\Field^n$. The unit entries of the vectors $\e_q\in\mathcal{Q}$ occupy the columns indexed by $\mathsf{S}'$, which are disjoint from the pivot columns of $\mathcal{N}$ indexed by $\mathsf{S}$. Using these unit entries as pivots eliminates the entries of the upper block $\mathcal{N}$ in the columns indexed by $\mathsf{S}'$. Since every $\e_q\in\mathcal{Q}$ vanishes outside its unit entry, these operations leave the $\mathsf{S}$-indexed identity columns invariant. The resulting matrix is therefore a permutation matrix and can be transformed into $\I_n$ by a suitable row permutation. Hence, the rows of $\widetilde{\R_{\mathsf{S}}}$ form a basis of $\Field^n$. Consequently, \[ \Field^n=\RS(\A)\oplus\Span(\mathcal{Q}), \] which yields \[ \Span(\mathcal{Q})=\Res_{\A}(\mathsf{S}'). \] \noindent\textit{ii)} Definition~\ref{def. RRPEF} implies that the $j$-th nonzero row of $\R_{\mathsf{S}}$ is the sum of its pivot unit vector $\e_{\mathcal{i}_j}$ and a linear combination of unit vectors indexed by $\mathsf{S}'$: \begin{equation}\label{eq. Noether basis formula} (\R_{\mathsf{S}})_{j\bullet}=r_{j\mathcal{i}_j}\e_{\mathcal{i}_j}+\sum_{q\in\mathsf{S}'}r_{jq}\e_q,\qquad j\in\I(\rho), \end{equation} where $\mathcal{i}_j\in\mathsf{S}$ and $r_{j\mathcal{i}_j}=1$. Equation~\eqref{eq. Noether basis formula} coincides with the Noether expansion~\eqref{Noether basis 1} for $\mathsf{C}=\mathsf{S}'$. Hence, $\mathcal{N}$ has the defining structure of the Noether basis associated with $\mathsf{C}$. By the uniqueness of the Noether basis $\Exnoether_{\mathsf{C}}$ relative to the CBI set $\mathsf{C}=\mathsf{S}'$, it follows that \[ \Exnoether_{\mathsf{C}}=\mathcal{N}, \] which completes the proof. \end{proof} Based on the equivalence relation defining $\mathcal{S}_{\A}^{o}$, we introduce the dual quotient set of CBI sets. For $\mathsf{C}_i,\mathsf{C}_j\in\mathcal{C}_{\A}$, with $\mathsf{S}_i=\mathsf{C}'_i$ and $\mathsf{S}_j=\mathsf{C}'_j$, define \[ \mathsf{C}_i\sim_{\mathcal{C}_{\A}}\mathsf{C}_j\quad\Longleftrightarrow\quad\R_{\mathsf{C}'_i}=\R_{\mathsf{C}'_j}. \] The corresponding quotient set is \[ \mathcal{C}_{\A}^{o}\triangleq\mathcal{C}_{\A}/\mathord{\sim_{\mathcal{C}_{\A}}}. \] The next theorem complements Proposition~\ref{prop. rows of RRPEF basis of RS} by constructing the RRPEF associated with an arbitrary Noether basis of $\RS(\A)$. Together, the two results establish the correspondence between the distinct Noether bases and the distinct RRPEFs. We first equip the Noether basis with the ordering required by Definition~\ref{def. RRPEF}. Let \[ k=\card(\Exnoether_{\mathsf{C}}). \] Since $\Exnoether_{\mathsf{C}}=\{\n_i\mid i\in\mathsf{C}'\}$ is a basis of $\RS(\A)$, one has \[ k=\card(\Exnoether_{\mathsf{C}})=\dim(\RS(\A))=\rank(\A)=\rho. \] Hence, \[ \card(\mathsf{C}')=\rho. \] Since $\mathsf{C}$ is a CBI set, the complementary criterion gives \[ \mathsf{S}=\mathsf{C}'\in\mathcal{S}_{\A}. \] Write the elements of $\mathsf{C}'$ in increasing order as \[ \mathsf{C}'=\{\beta_1<\cdots<\beta_{\rho}\}. \] Then \[ \beta:(\I(\rho),<)\longrightarrow(\mathsf{C}',<),\qquad j\longmapsto\beta_j, \] is the unique order isomorphism between these finite ordered sets. Applying this ordering to Eq.~\eqref{Noether basis 1}, with $\mathsf{S}=\mathsf{C}'$, gives \begin{equation}\label{eq. Noether formula2} (\Exnoether_{\mathsf{C}})_{j\bullet}=r_{j\beta_j}\e_{\beta_j}+\sum_{q\in\mathsf{C}}r_{jq}\e_q,\qquad j\in\I(\rho), \end{equation} where $\beta_j\in\mathsf{C}'$ and $r_{j\beta_j}=1$. Equation~\eqref{eq. Noether formula2} equips \[ \Exnoether_{\mathsf{C}}=\left\{(\Exnoether_{\mathsf{C}})_{j\bullet}\mid j\in\I(\rho)\right\} \] with the ordering induced by $\beta$. We refer to this indexed basis as the \emph{ordered Noether basis under $<$}. \begin{theorem}\label{theo. normalized RRPEF} Let $\A\in\Field^{m\times n}$ and let $\mathsf{C}\in\mathcal{C}_{\A}$. Let $\Exnoether_{\mathsf{C}}=\{\n_i\mid i\in\mathsf{C}'\}$ be the ordered Noether basis under $<$ for $\RS(\A)$. Construct the matrix obtained by stacking the $\rho=\rank(\A)$ elements of $\Exnoether_{\mathsf{C}}$ above $m-\rho$ zero rows: \begin{equation}\label{eq. RRPEF_stacking} \V_{\mathsf{C}'}\triangleq\left[\begin{array}{c}\Exnoether_{\mathsf{C}}\vspace{0.05in}\\ \hline\vspace{-0.1in}\\ \0\end{array}\right], \end{equation} where $\0$ denotes the $(m-\rho)\times n$ zero block. Then: \begin{enumerate}[\itshape i)] \item The matrix $\V_{\mathsf{C}'}$ is row equivalent to $\A$. \item The matrix $\V_{\mathsf{C}'}$ is the unique RRPEF of $\A$ relative to $\mathsf{S}=\mathsf{C}'$. Consequently, the set of distinct Noether bases \[ \left\{\Exnoether_{\mathsf{C}}\mid\mathsf{C}\in\mathcal{C}_{\A}\right\} \] is in bijection with $\Im(\mathcal{R}_{\A})$. \item The quotient sets $\mathcal{C}_{\A}^{o}$ and $\mathcal{S}_{\A}^{o}$ are in bijection. \end{enumerate} \end{theorem} \begin{proof} By the block construction of $\V_{\mathsf{C}'}$, \[ \nul_{\ell}(\V_{\mathsf{C}'})=m-\rho=\nul_{\ell}(\A), \] so $\V_{\mathsf{C}'}$ has exactly the number of zero rows required by the row nullity of $\A$. \noindent\textit{i)} Since $\Exnoether_{\mathsf{C}}$ is a basis of $\RS(\A)$, \[ \RS(\A)=\Span(\Exnoether_{\mathsf{C}}). \] The zero block in \eqref{eq. RRPEF_stacking} does not alter the row space, and therefore \[ \RS(\V_{\mathsf{C}'})=\Span(\Exnoether_{\mathsf{C}})=\RS(\A). \] Proposition~\ref{prop. equivalent statements to row equivalence} then gives \[ \V_{\mathsf{C}'}\REQ\A. \] \noindent\textit{ii)} By Eq.~\eqref{eq. Noether formula2}, $r_{j\beta_j}=1$ for every $j\in\I(\rho)$; hence $\V_{\mathsf{C}'}$ satisfies Postulate~(i) of Definition~\ref{def. RRPEF}. Postulate~(ii) follows from the order isomorphism $\beta$, while Postulate~(iii) follows from the block form in \eqref{eq. RRPEF_stacking}, whose lower block consists of exactly $\nul_{\ell}(\A)=m-\rho$ zero rows. For Postulate~(iv), Eq.~\eqref{eq. Noether formula2} shows that, in the $q$-th nonzero row, the only entry whose column index belongs to $\mathsf{C}'$ is the pivot entry $r_{q\beta_q}=1$. Hence, for $j\neq q$, \[ r_{q\beta_j}=0. \] Thus, each column indexed by $\beta_j\in\mathsf{C}'$ contains the unique nonzero entry $r_{j\beta_j}=1$, proving Postulate~(iv). Therefore, $\V_{\mathsf{C}'}$ is an RRPEF of $\A$ relative to \[ \mathsf{S}=\mathsf{C}'. \] By Theorem~\ref{theo.RNCSBT}, the Noether basis $\Exnoether_{\mathsf{C}}$ is uniquely determined by the CBI set $\mathsf{C}$. Its ordering is uniquely determined by the increasing order on $\mathsf{C}'$, and the lower zero block is fixed by $\nul_{\ell}(\A)=m-\rho$. Hence, $\V_{\mathsf{C}'}$ is the unique RRPEF of $\A$ relative to $\mathsf{S}=\mathsf{C}'$. Consequently, \[ \V_{\mathsf{C}'}=\R_{\mathsf{S}},\qquad\mathsf{S}=\mathsf{C}'. \] \noindent\textit{iii)} The complement operator induces the mapping \[ \Phi:\mathcal{C}_{\A}^{o}\longrightarrow\mathcal{S}_{\A}^{o},\qquad\Phi([\mathsf{C}])=[\mathsf{C}']. \] This mapping is well defined because, by the definitions of the two equivalence relations, \[ \mathsf{C}_i\sim_{\mathcal{C}_{\A}}\mathsf{C}_j\quad\Longleftrightarrow\quad\R_{\mathsf{C}'_i}=\R_{\mathsf{C}'_j}\quad\Longleftrightarrow\quad\mathsf{C}'_i\sim_{\mathcal{S}_{\A}}\mathsf{C}'_j. \] Since complementation is involutive, $(\mathsf{C}')'=\mathsf{C}$, the inverse mapping is \[ \Phi^{-1}:\mathcal{S}_{\A}^{o}\longrightarrow\mathcal{C}_{\A}^{o},\qquad\Phi^{-1}([\mathsf{S}])=[\mathsf{S}']. \] Therefore, $\Phi$ is bijective, which completes the proof. \end{proof} Proposition~\ref{prop. rows of RRPEF basis of RS} provides the inverse construction by recovering the Noether basis from the nonzero rows of the RRPEF. Consequently, the correspondence between distinct Noether bases and distinct RRPEFs is bijective. The preceding results are summarized in the following diagram: \begin{figure}[H] \centering \begin{tikzpicture}[node distance=2cm,auto] \node (C) {$\mathcal{C}_{\A}^{o}$}; \node (D) [right=5.5cm of C] {$\{\Exnoether_{\mathsf{C}}\mid\mathsf{C}\in\mathcal{C}_{\A}\}$}; \node (X) [below=1.5cm of C] {$\mathcal{S}_{\A}^{o}$}; \node (B) [below=1.5cm of D] {$\Im(\mathcal{R}_{\A})$}; \draw[->] (C) to node {Noether-basis bijection} (D); \draw[->] (C) to node [left] {Complement operator} (X); \draw[->] (D) to node [right=1ex] {Identifying bijection} (B); \draw[->] (X) to node {Quotient-induced bijection} (B); \end{tikzpicture} \caption{Commutative diagram} \label{fig:Commutative Arrow Diagram} \end{figure} The bijective correspondence between the Noether bases and the RRPEFs provides the structural characterization~\eqref{Noether basis 1} independently of a prescribed sequence of elementary row operations. Because the standard and Noether bases admit compatible extensions to row-finite matrices, the stacked construction in \eqref{eq. RRPEF_stacking} extends the RRPEF framework to row-finite settings. The corresponding algebraic framework is developed in~\cite{Pa:LRREFROWFINITE}.
 		
\subsection{Completeness of the RRPEF Family} The Rank Minor Theorem guarantees that a matrix of rank $\rho$ contains at least one nonsingular $\rho\times\rho$ minor. In a normalized row-equivalent form, the $\rho$ nonzero rows constitute a full-row-rank block, and the BCI collection $\mathcal{S}_{\A}$ indexes all nonsingular $\rho\times\rho$ transition submatrices of this block. The Rank Minor Theorem alone, however, does not determine a normalized coordinate representation. Such normalization is obtained by reducing each transition submatrix $\M$ to $\I_{\rho}$, thereby isolating the pivot coordinates prescribed by the corresponding BCI set $\mathsf{S}$. Thus, the Rank Minor Theorem provides the nonsingular algebraic foundation, whereas coordinate isolation yields the uniquely normalized representative $\R_{\mathsf{S}}$. We refer to this conjunction as \emph{full-rank coordinate isolation}. By Theorem~\ref{theo. normalized RRPEF}, each $\mathsf{S}\in\mathcal{S}_{\A}$ determines a unique RRPEF of $\A$. The following theorem establishes that these representatives exhaust all RRPEFs in the row-equivalence class of $\A$ and hence proves the completeness of $\mathcal{R}_{\A}$. \begin{theorem}[Completeness of the RRPEF Family]\label{theo. completeness_RRPEF} The indexed family \[ \mathcal{R}_{\A}=\bigl(\R_{\mathsf{S}}\bigr)_{\mathsf{S}\in\mathcal{S}_{\A}} \] is the complete canonical family of all \emph{RRPEFs} row-equivalent to $\A$. Specifically, if $\R$ is a row-normal form of $\A$ satisfying Definition~\ref{def. RRPEF}, then there exists an admissible index set $\mathsf{S}\in\mathcal{S}_{\A}$ such that \[ \R=\R_{\mathsf{S}}\in\mathcal{R}_{\A}. \] Equivalently, $\mathcal{R}_{\A}$ exhausts all \emph{RRPEF} representatives of the row-equivalence class of $\A$. \end{theorem} \begin{proof} Let $\R\in\Field^{m\times n}$ be an arbitrary RRPEF in the row-equivalence class of $\A$, and let $\C\in\Field^{\rho\times n}$ denote its full-row-rank nonzero block. By Definition~\ref{def. RRPEF}, $\C$ contains $\rho$ pivot columns forming the identity matrix. Let \[ \mathsf{S}=\{\mathcal{i}_1,\mathcal{i}_2,\ldots,\mathcal{i}_{\rho}\}\subseteq\I(n) \] be their column-index set. Then \[ \C_{\mathsf{S}}=\I_{\rho}. \] Let $\Ub\in\Field^{\rho\times n}$ denote the full-row-rank nonzero block of $\URREF(\A)$. Since $\R\REQ\A$, the matrices $\C$ and $\Ub$ have the same row space. Hence, there exists a nonsingular matrix $\Qb\in\Field^{\rho\times\rho}$ such that \[ \C=\Qb\Ub. \] Restricting this equality to the columns indexed by $\mathsf{S}$ gives \[ \I_{\rho}=\C_{\mathsf{S}}=\Qb\Ub_{\mathsf{S}}. \] Hence, \[ \Ub_{\mathsf{S}}=\Qb^{-1}, \] so $\Ub_{\mathsf{S}}$ is nonsingular. By the Column-Index Pivot Criterion, \[ \mathsf{S}\in\mathcal{S}_{\A}. \] Moreover, \[ \Qb=\Ub_{\mathsf{S}}^{-1}, \] and consequently, \[ \C=\Qb\Ub=\Ub_{\mathsf{S}}^{-1}\Ub. \] Thus, $\C$ is the uniquely normalized full-row-rank block associated with the BCI set $\mathsf{S}$. Appending the prescribed $m-\rho$ zero rows yields \[ \R=\R_{\mathsf{S}}. \] Therefore, every RRPEF in the row-equivalence class of $\A$ is indexed by an admissible BCI set $\mathsf{S}\in\mathcal{S}_{\A}$ and belongs to $\mathcal{R}_{\A}$. Hence, $\mathcal{R}_{\A}$ is complete. \end{proof} The RRPEF of $\A\in\Field^{m\times n}$ relative to a BCI set $\mathsf{S}$ is constructed by reducing the transition submatrix $\Ub_{\mathsf{S}}$ of the full-row-rank nonzero block of $\URREF(\A)$ to $\I_{\rho}$ and appending the prescribed $m-\rho$ zero rows. The resulting matrix is row equivalent to $\A$, has the prescribed pivot columns, and contains exactly $m-\rho$ zero rows at its lower boundary. Constraint matrices in simplex tableaux provide standard examples of this construction.

\subsection{A Worked Example of RRPEFs}\label{subsec: Worked Example RRPEF} The following example demonstrates the effect of selecting different CBI sets on the formation of RRPEFs and illustrates how the complementary criterion identifies admissible pivot sets. \begin{example}\label{ex. WorkedExampleForRRPEFs} \normalfont Consider the matrix \begin{equation}\label{eq. 5x6 matrix} \A=\left[\begin{array}{cccccc} 1&0&a&0&0&h\\ 0&1&b&0&0&d\\ 0&0&0&1&0&f\\ 0&0&0&0&1&e\\ 0&0&0&0&0&0 \end{array}\right], \end{equation} where $a,b,h,d,e,f$ are nonzero and \[ \det\begin{bmatrix}a&h\\b&d\end{bmatrix}\neq0. \] The case \[ \det\begin{bmatrix}a&h\\b&d\end{bmatrix}=0 \] is presented in Appendix~\ref{append. MatrixSingularity}. \medskip The matrix $\A$ is in URREF and \[ \A=\R_{\{1,2,4,5\}}. \] Thus, its BCI set is \[ \mathsf{S}=\{1,2,4,5\}, \] and the associated CBI set is \[ \mathsf{C}=\mathsf{S}'=\{3,6\}. \] Because $\rank(\A)=4$, the residual space $\Res_{\A}(\{3,6\})$ has dimension $2$, equal to the codimension of $\RS(\A)$.\\
 \medskip 
 \noindent\textbf{Step 1: Complementary criterion via the augmented determinant test.} To obtain a non-extreme RRPEF, consider the candidate CBI set $\mathsf{C}=\{2,6\}$, whose complementary BCI set is $\mathsf{S}=\{1,3,4,5\}$. To verify that $\{\e_2,\e_6\}$ spans a complementary subspace of $\RS(\A)$, we append these vectors to the four-row full-rank nonzero block of $\A$; see the proof of Proposition~\ref{prop. rows of RRPEF basis of RS}: \begin{equation}\label{eq. augmented 5x6 matrix} \boldsymbol{\rm T}=\left[\begin{array}{cccccc} 1&0&a&0&0&h\\ 0&1&b&0&0&d\\ 0&0&0&1&0&f\\ 0&0&0&0&1&e\\ 0&1&0&0&0&0\\ 0&0&0&0&0&1 \end{array}\right]. \end{equation} Since \[ \det(\boldsymbol{\rm T})=-b\neq0, \] the rows of $\boldsymbol{\rm T}$ form a basis of $\Field^6$. Hence, $\mathsf{C}=\{2,6\}$ is an admissible CBI set. 
	
\medskip 

\noindent\textbf{Step 2: Constructing non-extreme RRPEFs.} Reducing the transition submatrix indexed by \[ \mathsf{S}=\{1,3,4,5\}\in\mathcal{S}_{\A} \] to $\I_4$ by Gauss--Jordan elimination yields the non-extreme RRPEF \[ \R_{\{1,3,4,5\}}=\begin{bmatrix} 1&-\frac{a}{b}&0&0&0&\frac{hb-ad}{b}\\ 0&\frac{1}{b}&1&0&0&\frac{d}{b}\\ 0&0&0&1&0&f\\ 0&0&0&0&1&e\\ 0&0&0&0&0&0 \end{bmatrix}. \] Similarly, augmenting the full-row-rank nonzero block of $\A$ with $\{\e_2,\e_5\}$ yields another admissible CBI set, \[ \mathsf{C}=\{2,5\}, \] whose complementary BCI set is \[ \mathsf{S}=\{1,3,4,6\}. \] Reducing the corresponding transition submatrix to $\I_4$ gives \begin{equation}\label{eq. matrix C} \R_{\{1,3,4,6\}}=\left[\begin{array}{cccccc} 1&-\frac{a}{b}&0&0&\frac{ad-hb}{eb}&0\\ 0&\frac{1}{b}&1&0&-\frac{d}{eb}&0\\ 0&0&0&1&-\frac{f}{e}&0\\ 0&0&0&0&\frac{1}{e}&1\\ 0&0&0&0&0&0 \end{array}\right]. \end{equation} \medskip \noindent\textbf{Step 3: Counting admissible RRPEFs.} Since \[ \binom{6}{2}=\binom{6}{4}=15, \] there are $15$ candidate CBI sets, equivalently, $15$ complementary candidate BCI sets of cardinality four. Only the CBI sets \[ \{5,6\},\qquad\{4,6\},\qquad\{4,5\} \] fail the complementary criterion. Their complements are therefore not BCI sets. Consequently, exactly $12$ of the $15$ candidate BCI sets are admissible and produce RRPEFs of $\A$. In this example, the admissible BCI sets generate distinct RRPEFs; hence, \[ \mathcal{S}_{\A}^{o}=\mathcal{S}_{\A}. \] \medskip \noindent\textbf{Step 4: Constructing the LRREF.} The maximal RRPEF, relative to the BCI set \[ \mathsf{S}=\{3,4,5,6\}, \] is \begin{equation}\label{equation LRREF of 5x6 matrix} \R_{\{3,4,5,6\}}=\left[\begin{array}{cccccc} \frac{d}{ad-bh}&-\frac{h}{ad-bh}&1&0&0&0\\ \frac{bf}{ad-bh}&-\frac{af}{ad-bh}&0&1&0&0\\ \frac{be}{ad-bh}&-\frac{ae}{ad-bh}&0&0&1&0\\ -\frac{b}{ad-bh}&\frac{a}{ad-bh}&0&0&0&1\\ 0&0&0&0&0&0 \end{array}\right], \end{equation} which is the LRREF of $\A$. \end{example}
 
\section{The Lexicographic Hierarchy of RRPEFs}\label{sec. The Hierarchy of normalized Forms} With respect to the lexicographic order introduced in Subsection~\ref{subsec: Lexicographic Order on BCI Sets}, the URREF and LRREF constitute the two extreme normalized forms, while the remaining RRPEFs are classified as \emph{internal forms}. This section develops the lexicographic structure of the RRPEF family and establishes the strict total order transferred from the core BCI sets to the distinct RRPEFs. \subsection{The Extreme RRPEFs} The leftmost and rightmost pivot sets arise as the two extreme admissible choices within the RRPEF family. They correspond, respectively, to the minimal and maximal elements of $\mathcal{S}_{\A}$ and determine the URREF and LRREF. \begin{definition}[Extreme pivot positions]\label{def. LPP and RPP} The element \[ \mathsf{S}_{\ell}=\{\mathcal{i}_1,\ldots,\mathcal{i}_{\rho}\}\in\mathcal{S}_{\A} \] is called a \emph{leftmost pivot position} \emph{(LPP)} set if, for every \[ \mathsf{S}=\{q_1,\ldots,q_{\rho}\}\in\mathcal{S}_{\A}\setminus\{\mathsf{S}_{\ell}\}, \] one has \[ \mathcal{i}_j\leq q_j\quad\text{for all }j\in\I(\rho),\qquad\mathcal{i}_i<q_i\quad\text{for at least one }i\in\I(\rho). \] Similarly, the element \[ \mathsf{S}_{\calr}=\{\mathcal{j}_1,\ldots,\mathcal{j}_{\rho}\}\in\mathcal{S}_{\A} \] is called a \emph{rightmost pivot position} \emph{(RPP)} set if, for every \[ \mathsf{S}=\{q_1,\ldots,q_{\rho}\}\in\mathcal{S}_{\A}\setminus\{\mathsf{S}_{\calr}\}, \] one has \[ \mathcal{j}_k\geq q_k\quad\text{for all }k\in\I(\rho),\qquad\mathcal{j}_i>q_i\quad\text{for at least one }i\in\I(\rho). \] \end{definition} An RRPEF relative to the LPP set $\mathsf{S}_{\ell}$, respectively the RPP set $\mathsf{S}_{\calr}$, is said to satisfy the leftmost, respectively rightmost, pivot criterion. \begin{lemma}[Lexicographic Minimality of the $\mathsf{S}_{\ell}$ Set]\label{lem:minimalityURREF} Let \[ \mathsf{S}_{\ell}=\{\mathcal{i}_1,\ldots,\mathcal{i}_{\rho}\}\in\mathcal{S}_{\A} \] be the admissible pivot-index set obtained by Gaussian elimination through the successive selection of the leftmost admissible pivot. Then \begin{equation}\label{eq. minimal Sell} \mathsf{S}_{\ell}=\min_{\prec_L}\mathcal{S}_{\A}. \end{equation} \end{lemma} \begin{proof} Suppose that there exists \[ \mathsf{S}=\{q_1,\ldots,q_{\rho}\}\in\mathcal{S}_{\A} \] such that \[ \mathsf{S}\prec_L\mathsf{S}_{\ell}. \] Let $h$ be the first index at which their increasing representations differ. Then \[ q_i=\mathcal{i}_i,\qquad i<h, \] and \[ q_h<\mathcal{i}_h. \] Since $\mathsf{S}$ is an admissible BCI set, its first $h$ columns extend the previously selected $h-1$ pivot columns to a linearly independent set. Consequently, after the first $h-1$ elimination stages, an admissible pivot exists in a column no farther to the right than $q_h$. This contradicts the selection of $\mathcal{i}_h$ as the leftmost admissible pivot at the $h$-th stage. Therefore, no lexicographically smaller admissible BCI set exists, proving~\eqref{eq. minimal Sell}. \end{proof} \begin{proposition}\label{prop. Extreme Properties of RREFs} The \emph{URREF} and \emph{LRREF} of $\A$ are, respectively, the minimal and maximal elements of $\Im(\mathcal{R}_{\A})$ under the transferred lexicographic order. \end{proposition} \begin{proof} By Lemma~\ref{lem:minimalityURREF}, the pivot set $\mathsf{S}_{\ell}$ generated by Gaussian elimination is the lexicographically least element of $\mathcal{S}_{\A}$. The RRPEF corresponding to $\mathsf{S}_{\ell}$ is the URREF of $\A$. Hence, \[ \R_{\mathrm{URREF}}=\min_{\prec_L}\Im(\mathcal{R}_{\A}). \] Dually, rightmost pivot selection yields the lexicographically greatest admissible pivot-index set $\mathsf{S}_{\calr}$, whose corresponding RRPEF is the LRREF. Consequently, \[ \R_{\mathrm{LRREF}}=\max_{\prec_L}\Im(\mathcal{R}_{\A}), \] which completes the proof. \end{proof} The LPP set is the lexicographically least BCI set and determines the URREF. Dually, the RPP set is the lexicographically greatest BCI set and determines the LRREF. The remaining minimal representatives of the equivalence classes determine the internal RRPEFs. The lexicographic order is defined on $\binom{[n]}{\rho}$ and restricts to the set of minimal representatives identified with $\mathcal{S}_{\A}^{o}$. The quotient-induced bijection \[ \mathcal{S}_{\A}^{o}\longleftrightarrow\Im(\mathcal{R}_{\A}) \] in Figure~\ref{fig:Commutative Arrow Diagram} then transfers this order to the set of distinct RRPEFs. The following theorem summarizes the preceding results and incorporates the transferred lexicographic order into the RRPEF framework. \begin{theorem}[The Lexicographic RRPEF Chain Theorem]\label{theo. The Lexicographic RRPEF Chain} The set $\Im(\mathcal{R}_{\A})$ of distinct \emph{RRPEFs} of $\A$ satisfies the following properties: \begin{enumerate}[\itshape i)] \item The bijection between $\mathcal{S}_{\A}^{o}$ and $\Im(\mathcal{R}_{\A})$ transfers the strict lexicographic total order of $(\mathcal{S}_{\A}^{o},\prec_L)$ to $\Im(\mathcal{R}_{\A})$: \[ \mathsf{S}_m\prec_L\mathsf{S}_{m+1}\quad\Longleftrightarrow\quad\R_{\mathsf{S}_m}\prec_{\mathcal{R}_{\A}}\R_{\mathsf{S}_{m+1}}. \] Since $\prec_{\mathcal{R}_{\A}}$ is transferred directly from $\prec_L$, the notation $\prec_L$ is used for both orders. \item The matrices $\R_{\mathrm{URREF}}$ and $\R_{\mathrm{LRREF}}$ belong to $\Im(\mathcal{R}_{\A})$ and constitute the unique extreme elements of the chain: \begin{equation}\label{eq. poles} \R_{\mathrm{URREF}}=\min_{\prec_L}\Im(\mathcal{R}_{\A})\qquad\text{and}\qquad\R_{\mathrm{LRREF}}=\max_{\prec_L}\Im(\mathcal{R}_{\A}). \end{equation} Moreover, for every $\R_{\mathsf{S}}\in\Im(\mathcal{R}_{\A})$, \[ \R_{\mathrm{URREF}}\preceq_L\R_{\mathsf{S}}\preceq_L\R_{\mathrm{LRREF}}. \] \end{enumerate} \end{theorem} \noindent In constructive terms, the admissible BCI sets of the matrix $\A$ in \eqref{eq. 5x6 matrix} are determined by nonsingular coordinate isolation. Gaussian elimination selects pivots at $(1,1)$ and $(2,2)$, the leftmost nonzero entries of the first two rows. In row~3, the entry in column~3 vanishes, so the pivot search proceeds to column~4 and selects $(3,4)$. Row~4 then yields the pivot $(4,5)$. Thus, the LPP set is \[ \mathsf{S}_{\ell}=\{1,2,4,5\}, \] which is the lexicographically least BCI set of $\A$. The three lexicographically smaller candidates in $\binom{[6]}{4}$ satisfy \[ \{1,2,3,4\}\prec_L\{1,2,3,5\}\prec_L\{1,2,3,6\}\prec_L\{1,2,4,5\} \] and are inadmissible. Indeed, if one of the first three candidates were admissible, it would contradict the leftmost-pivot construction of $\mathsf{S}_{\ell}=\{1,2,4,5\}$. The LPP set is therefore the first admissible element of $\mathcal{S}_{\A}$, whereas the RPP set is its last admissible element. Alternatively, using the relation $\mathsf{S}=\mathsf{C}'$, the BCI chain may be constructed indirectly from the collection of CBI sets by applying the augmented-matrix determinant test to the candidates in $\binom{[n]}{n-\rho}$ in reverse lexicographic order; see Appendix~\ref{append.AdmissibilityOfBCIViaCBISets}. The matrix $\A$ admits the following $12$ BCI sets: \begin{align*} \{1,2,4,5\}&\prec_L\{1,2,4,6\}\prec_L\{1,2,5,6\}\prec_L\{1,3,4,5\}\prec_L\{1,3,4,6\}\prec_L\{1,3,5,6\}\\ &\prec_L\{1,4,5,6\}\prec_L\{2,3,4,5\}\prec_L\{2,3,4,6\}\prec_L\{2,3,5,6\}\prec_L\{2,4,5,6\}\prec_L\{3,4,5,6\}. \end{align*} These sets yield the corresponding $12$ distinct RRPEFs. In this example, \[ \mathcal{S}_{\A}^{o}=\mathcal{S}_{\A}. \] \subsection{The RRPEF Chain and Its Extreme Elements}\label{subsec:RRPEFChain} The extreme elements of the RRPEF chain are obtained constructively within the Gauss--Jordan reduction framework. Let $\Per_n$ be the reversal permutation matrix of order $n$, and let $\Per_{k,m}$ denote the partial reversal permutation matrix \begin{align*} \Per_n&\triangleq\left[\begin{array}{cccc}0&0&\cdots&1\\[2pt]\vdots&\vdots&\iddots&\vdots\\[6pt]1&0&\cdots&0\end{array}\right],\qquad \Per_{k,m}\triangleq\left[\begin{array}{cc}\Per_k&\0\\[6pt]\0&\I_{m-k}\end{array}\right]. \end{align*} In particular, $\Per_{\rho,m}$ reverses the first $\rho$ rows while leaving the bottom $m-\rho$ rows fixed. Applying this operation to $\URREF(\A\Per_n)$ and then restoring the original column order produces $\LRREF(\A)$. \begin{proposition}[Construction of the terminal pole]\label{prop:LRREFconstruction} Let $\A\in\Field^{m\times n}$ have rank $\rho$. Then \begin{equation}\label{eq.LRREFconstruction} \LRREF(\A)=\Per_{\rho,m}\Bigl(\URREF(\A\Per_n)\Bigr)\Per_n. \end{equation} \end{proposition} \begin{proof} Let \[ \Ub=\URREF(\A),\qquad\Lb=\LRREF(\A),\qquad\text{and}\qquad\Vb=\URREF(\A\Per_n). \] Define \begin{equation}\label{eq.LRREFconstruction2} \Lb^{*}=\Per_{\rho,m}\Vb\Per_n. \end{equation} Since $\Vb$ is in URREF, reversing its columns and then reversing its $\rho$ nonzero rows transforms the leftmost-pivot structure into the rightmost-pivot structure while preserving the zero rows at the bottom. Hence, $\Lb^{*}$ is in LRREF. It remains to show that $\Lb^{*}\REQ\A$. Since $\Vb\REQ\A\Per_n$, there exists an invertible matrix $\N$ such that \begin{equation}\label{equation U and AP row-equivalent} \N(\A\Per_n)=\Vb. \end{equation} Because \[ \Per_n=\Per_n^{\top}=\Per_n^{-1}\qquad\text{and}\qquad\Per_{\rho,m}^{-1}=\Per_{\rho,m}, \] Eq.~\eqref{eq.LRREFconstruction2} gives \[ \Vb=\Per_{\rho,m}\Lb^{*}\Per_n. \] Substituting this expression into \eqref{equation U and AP row-equivalent} and post-multiplying by $\Per_n$ yields \[ \N\A\Per_n=\Per_{\rho,m}\Lb^{*}\Per_n\quad\Longleftrightarrow\quad\N\A=\Per_{\rho,m}\Lb^{*}. \] Pre-multiplication by $\Per_{\rho,m}$ gives \[ (\Per_{\rho,m}\N)\A=\Lb^{*}. \] Since \[ \Qb\triangleq\Per_{\rho,m}\N \] is invertible, $\Lb^{*}\REQ\A$. The matrix $\Lb^{*}$ is therefore an LRREF of $\A$, and the uniqueness of the LRREF implies \[ \Lb^{*}=\Lb, \] as claimed. \end{proof} Consequently, the terminal element of the RRPEF chain is obtained from the same Gauss--Jordan procedure as the initial element, after reversing the column order and the order of the nonzero rows. Identity~\eqref{eq.LRREFconstruction} establishes the constructive duality between $\URREF(\A)$ and $\LRREF(\A)$ and confirms their status as the two row-canonical extreme members of the complete RRPEF family. \subsection{Cardinality Bounds for the RRPEF Chain} Let $\A\in\Field^{m\times n}$ have rank $\rho$. The cardinality $\card(\mathcal{S}_{\A})$ counts the admissible BCI sets, whereas $\card(\mathcal{S}_{\A}^{o})$ counts the distinct RRPEFs. Thus, \[ \card(\mathcal{S}_{\A}^{o})=\card\bigl(\Im(\mathcal{R}_{\A})\bigr)\leq\card(\mathcal{S}_{\A}). \] If $\rho=m=n$, then $\A$ is nonsingular and \[ \card(\mathcal{S}_{\A})=\card(\mathcal{S}_{\A}^{o})=1. \] If $\rho=n\leq m$, then $\I(n)$ is the unique $n$-element subset of itself, and every RRPEF has $\I_n$ as its nonzero block. Hence, \[ \card(\mathcal{S}_{\A})=\card(\mathcal{S}_{\A}^{o})=1. \] For $0<\rho<n$, multiple admissible pivot-index sets may occur, with the URREF and LRREF providing the two extremes. In general, \[ 1\leq\card(\mathcal{S}_{\A}^{o})\leq\card(\mathcal{S}_{\A})\leq\binom{n}{\rho}. \] \section{The Induced Lexicographic Canonical Algorithm}\label{sec. ILCA} In this section, we introduce the \emph{Induced Lexicographic Canonical Algorithm} (ILCA), which constructs the finite lexicographic chain of the RRPEF family. The two extreme RRPEFs are obtained automatically through the leftmost- and rightmost-pivot conventions. The internal forms are generated according to the strict lexicographic ranking transferred from the core BCI sets $\mathcal{S}_{\A}^{o}$. Thus, ILCA provides an algebraic traversal from the URREF to the LRREF through the distinct internal RRPEFs. The traversal is implemented through the order isomorphism \[ \binom{[n]}{\rho}\cong\I\left(\binom{n}{\rho}\right), \] which identifies each candidate pivot-index set with its position in the strict lexicographic order. \subsection{Lexicographic Ranking and Successor Generation} The \emph{lexicographic rank function} \begin{equation} \label{eq. lexicographic rank} f:\binom{[n]}{\rho}\longrightarrow\I\left(\binom{n}{\rho}\right),\qquad f(\mathcal{i}_1,\ldots,\mathcal{i}_{\rho}) = 1+\sum_{j=1}^{\rho}\sum_{t=\mathcal{i}_{j-1}+1}^{\mathcal{i}_j-1}\binom{n-t}{\rho-j}, \qquad \mathcal{i}_0\triangleq0, \end{equation} enumerates $\binom{[n]}{\rho}$ in strict lexicographic order and provides an order isomorphism between \[ \left(\binom{[n]}{\rho},\prec_L\right)\qquad\text{and}\qquad\left(\I\left(\binom{n}{\rho}\right),<\right) \] \cite{BradManSip2005}. For $\mathsf{S}\in\binom{[n]}{\rho}$ with \[ f(\mathsf{S})<\binom{n}{\rho}, \] the unranking map $f^{-1}$ gives \[ f(\mathsf{S})<f(\mathsf{S})+1\quad\Longleftrightarrow\quad f^{-1}\bigl(f(\mathsf{S})\bigr)\prec_Lf^{-1}\bigl(f(\mathsf{S})+1\bigr)\quad\Longleftrightarrow\quad\mathsf{S}\prec_L\mathsf{S}_{\mathrm{next}}. \] Accordingly, the immediate lexicographic successor $\mathsf{S}_{\mathrm{next}}$ of $\mathsf{S}$ is \begin{equation}\label{eq. Snext} \mathsf{S}_{\mathrm{next}}=f^{-1}\bigl(f(\mathsf{S})+1\bigr). \end{equation} The inverse map $f^{-1}$ can be evaluated by standard lexicographic unranking algorithms for combinations; see~\cite{GenitriniPepin2021}. \subsection{Transition-Contraction Mechanisms} 
\begin{figure}[htbp] \begin{center} \begin{tikzpicture}[ scale=1.0, transform shape, node distance=4.5mm, every node/.style={font=\small}, startstop/.style={rectangle,rounded corners,draw=blue!70,fill=blue!15,minimum width=32mm,minimum height=8mm,align=center,line width=0.8pt}, node distance=0.6cm and 1.5cm, every node/.style={font=\small,inner sep=4pt}, process/.style={rectangle,minimum width=2.2cm,minimum height=0.6cm,align=center,draw=black!70,fill=blue!10,rounded corners}, decision/.style={diamond,minimum width=2.2cm,minimum height=1.2cm,align=center,draw=black!70,fill=orange!10,aspect=1.6}, startstop/.style={ellipse,minimum width=2.0cm,minimum height=0.6cm,align=center,draw=black!70,fill=green!15}, io/.style={trapezium,trapezium left angle=70,trapezium right angle=110,minimum width=2.2cm,minimum height=0.6cm,align=center,draw=black!70,fill=purple!10}, preparation/.style={chamfered rectangle,chamfered rectangle sep=4pt,minimum width=2.2cm,minimum height=0.6cm,align=center,draw=black!70,fill=yellow!15}, arrow/.style={thick,-{Stealth}} ] \node (startnode) [startstop] {START}; \node (readA) [io,below=of startnode] {Read $\A$}; \node (start) [preparation,below=of readA] {Calculate the $\URREF(\A)$\\Set $\R_{\mathsf{S}_{\ell}}\equiv\URREF(\A)$;\\Initial values: $\mathcal{S}_{\A}\leftarrow\{\mathsf{S}_{\ell}\}$; $\mathcal{R}_{\A}\leftarrow\{(\R_{\mathsf{S}_{\ell}},\mathsf{S}_{\ell})\}$; $\mathsf{S}\leftarrow\mathsf{S}_{\ell}$}; \node (step2) [process,below=0.7cm of start] {Identify the immediate successor\\[0.01in]$\mathsf{S}_{\mathrm{next}}$ of $\mathsf{S}$ within $\binom{[n]}{\rho}$ relative to $\prec_L$\\[0.01in]Set $\mathsf{S}\leftarrow\mathsf{S}_{\mathrm{next}}$}; \node (step3) [process,below=of step2] {Isolate the $\rho\times\rho$ submatrix $\M$\\of $\R_{\mathsf{S}_{\ell}}$ with column indices in $\mathsf{S}$}; \node (decide) [decision,below=of step3] {Is $\M$ invertible?\\$\det(\M)\neq0$}; \node (assign) [process,below=of decide] {Set $\mathcal{S}_{\A}\leftarrow\mathcal{S}_{\A}\cup\{\mathsf{S}\}$}; \node (step4) [process,below=of assign] {Transform $\R_{\mathsf{S}_{\ell}}$ into $\R_{\mathsf{S}}$ by\\the reduction of the nonsingular\\transition submatrix $\M$ to the identity $\I_{\rho}$\\$\mathcal{R}_{\A}\leftarrow\mathcal{R}_{\A}\cup\{(\R_{\mathsf{S}},\mathsf{S})\}$}; \node (step5) [decision,below=of step4] {Has the chain reached\\the terminal envelope?\\$\mathsf{S}=\mathsf{S}_{\calr}$?}; \node (endtext) [io,below=of step5] {Output: The complete indexed family $\mathcal{R}_{\A}$ of RRPEFs}; \node (endnode) [startstop,below=of endtext] {END}; \draw[arrow] (startnode) -- (readA); \draw[arrow] (readA) -- (start); \draw[arrow] (start) -- (step2); \draw[arrow] (step2) -- (step3); \draw[arrow] (step3) -- (decide); \draw[arrow] (decide.180) -- node[above] {\scriptsize NO} ++(-1.5,0) |- (step2.180); \draw[arrow] (decide.270) -- node[right] {\scriptsize YES} (assign); \draw[arrow] (assign) -- (step4); \draw[arrow] (step4) -- (step5); \draw[arrow] (step5.0) -- node[above] {\scriptsize NO} ++(1.8,0) |- (step2.0); \draw[arrow] (step5.270) -- node[right] {\scriptsize YES} (endtext); \draw[arrow] (endtext) -- (endnode); \end{tikzpicture} \caption{ILCA Trajectory} \label{fig:Algorithm}\vspace{-0.1in} \end{center} \end{figure}
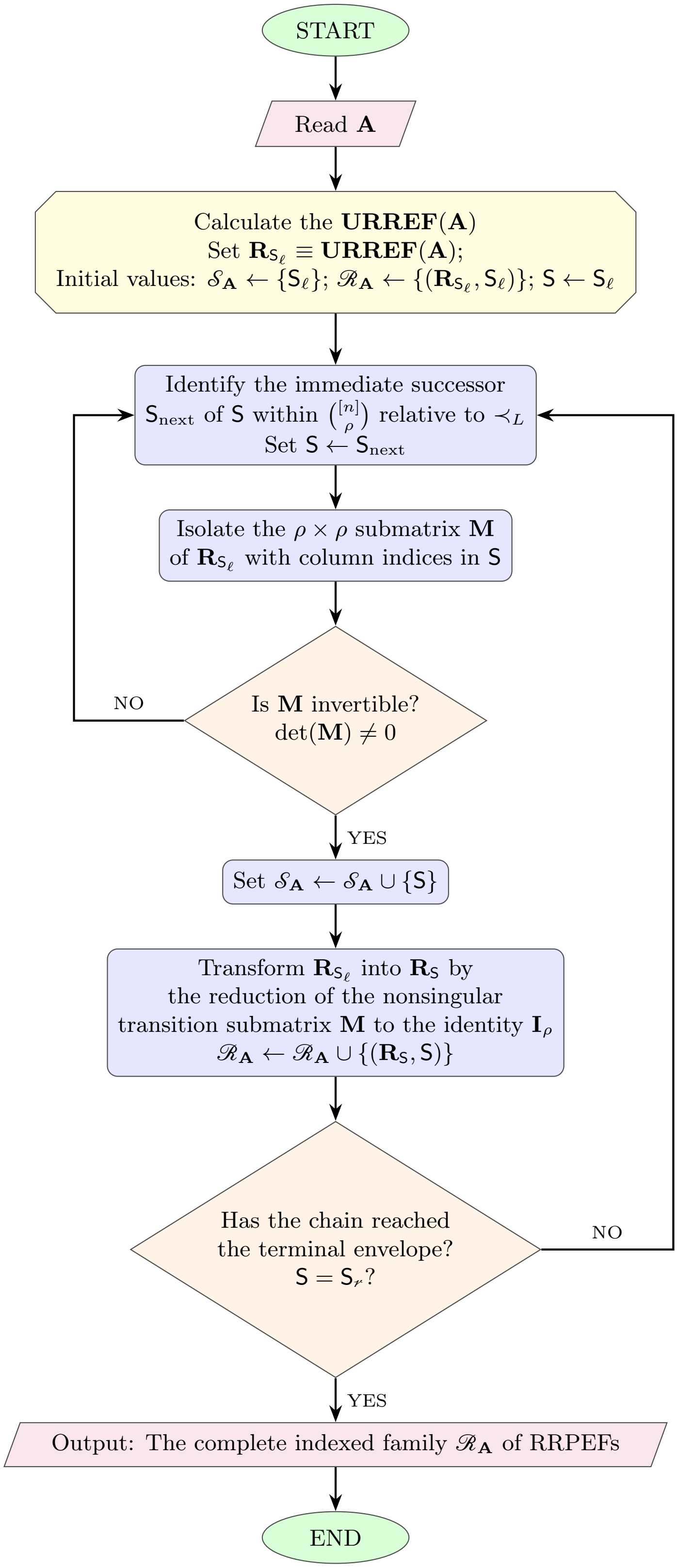

The following proposition establishes the transition-contraction mechanisms underlying the ILCA traversal: lexicographic confinement, singular-transition skipping, duplicate-RRPEF reuse, and scalar-equivalent RRPEF recovery. \begin{proposition}[Transition-contraction mechanisms]\label{prop. transition reuse} Let $\A\in\Field^{\rho\times n}$ have full row rank, and let $\mathsf{S}_i,\mathsf{S}_j\in\binom{[n]}{\rho}$ determine the transition submatrices \[ \M_i=\A_{\mathsf{S}_i},\qquad\M_j=\A_{\mathsf{S}_j}. \] The following statements hold: \begin{enumerate}[\itshape i)] \item \emph{Lexicographic confinement.} Every candidate pivot-index set lying lexicographically outside the interval \[ [\mathsf{S}_{\ell},\mathsf{S}_{\calr}]_L \] determines a singular transition submatrix and is excluded from the \emph{ILCA} traversal. \item \emph{Singular-transition skipping.} If \[ \det(\M_j)=0, \] then $\mathsf{S}_j\notin\mathcal{S}_{\A}$; hence, no \emph{RRPEF} relative to $\mathsf{S}_j$ exists, and the candidate $\mathsf{S}_j$ is skipped. \item \emph{Duplicate-RRPEF reuse.} If \[ \M_j=\M_i, \] then \[ \R_{\mathsf{S}_j}=\R_{\mathsf{S}_i}. \] Consequently, the previously constructed \emph{RRPEF} is reused, and no additional matrix inversion or Gauss--Jordan reduction is required. \item \emph{Scalar-equivalent RRPEF recovery.} If \[ \M_j=\lambda\M_i \] for some $\lambda\in\Field\setminus\{0,1\}$, then \[ \R_{\mathsf{S}_j}=\lambda^{-1}\R_{\mathsf{S}_i}. \] Although $\R_{\mathsf{S}_j}$ is generally distinct from $\R_{\mathsf{S}_i}$, it is recovered by direct scalar rescaling, without a new matrix inversion or Gauss--Jordan reduction. \end{enumerate} \end{proposition} \begin{proof} \noindent\textit{i)} By Lemma~\ref{lem:minimalityURREF} and its dual counterpart for the LRREF, \[ \mathsf{S}_{\ell}=\min_{\prec_L}\mathcal{S}_{\A}\qquad\text{and}\qquad\mathsf{S}_{\calr}=\max_{\prec_L}\mathcal{S}_{\A}. \] Let $\mathsf{S}\in\binom{[n]}{\rho}$ satisfy \[ \mathsf{S}\prec_L\mathsf{S}_{\ell}\qquad\text{or}\qquad\mathsf{S}_{\calr}\prec_L\mathsf{S}. \] If $\A_{\mathsf{S}}$ were nonsingular, then $\mathsf{S}\in\mathcal{S}_{\A}$, contradicting the minimality of $\mathsf{S}_{\ell}$ or the maximality of $\mathsf{S}_{\calr}$. Hence, \[ \det(\A_{\mathsf{S}})=0. \] Therefore, every candidate outside the indicated interval is excluded by lexicographic confinement. \noindent\textit{ii)} By Definition~\ref{def. basis_indexing_set}, \[ \mathsf{S}_j\in\mathcal{S}_{\A}\quad\Longleftrightarrow\quad\det(\M_j)\neq0. \] Therefore, if $\det(\M_j)=0$, then $\mathsf{S}_j\notin\mathcal{S}_{\A}$. Since an RRPEF is defined only relative to an admissible BCI set, no RRPEF relative to $\mathsf{S}_j$ exists, and the candidate is removed by singular-transition skipping. \noindent\textit{iii)} Suppose that $\M_i$ and $\M_j$ are nonsingular. The corresponding RRPEFs are \[ \R_{\mathsf{S}_k}=\M_k^{-1}\A,\qquad k\in\{i,j\}. \] If $\M_j=\M_i$, then $\M_j^{-1}=\M_i^{-1}$, and therefore \[ \R_{\mathsf{S}_j}=\M_j^{-1}\A=\M_i^{-1}\A=\R_{\mathsf{S}_i}. \] Hence, the existing RRPEF is reused, and no additional construction is required. \noindent\textit{iv)} If $\M_j=\lambda\M_i$, where $\lambda\in\Field\setminus\{0,1\}$, then \[ \M_j^{-1}=(\lambda\M_i)^{-1}=\lambda^{-1}\M_i^{-1}. \] Consequently, \[ \R_{\mathsf{S}_j}=\M_j^{-1}\A=\lambda^{-1}\M_i^{-1}\A=\lambda^{-1}\R_{\mathsf{S}_i}. \] Thus, $\R_{\mathsf{S}_j}$ is recovered by scalar rescaling, which completes the proof. \end{proof} \noindent Here, transition contraction refers to the reduction of explicit RRPEF constructions. Lexicographic confinement restricts the candidate domain, singular-transition skipping removes candidates that generate no RRPEF, duplicate-RRPEF reuse avoids reconstructing an existing RRPEF, and scalar-equivalent RRPEF recovery obtains a distinct RRPEF by direct rescaling. The indexed RRPEF family admits the set-theoretic realization \[ \mathcal{R}_{\A}\triangleq\left\{(\R_{\mathsf{S}},\mathsf{S})\mid\mathsf{S}\in\mathcal{S}_{\A}\right\}. \] Thus, distinct BCI sets remain distinct indices of the family even when they generate the same RRPEF matrix. The corresponding projection onto the image set is \[ \pi:\mathcal{R}_{\A}\longrightarrow\Im(\mathcal{R}_{\A}),\qquad\pi(\R_{\mathsf{S}},\mathsf{S})=\R_{\mathsf{S}}. \] The complete indexed RRPEF family can be generated by the basic ILCA trajectory displayed in Figure~\ref{fig:Algorithm}. This streamlined version explicitly implements lexicographic confinement and singular-transition skipping: the traversal begins at $\mathsf{S}_{\ell}$, terminates at $\mathsf{S}_{\calr}$, and bypasses every candidate whose transition submatrix is singular. For visual clarity, duplicate-RRPEF reuse and scalar-equivalent RRPEF recovery are not represented in the diagram. Accordingly, the displayed version may reconstruct the same RRPEF matrix from more than one BCI set while retaining the corresponding indexed pairs in $\mathcal{R}_{\A}$. At the level of distinct matrices, the full algebraic ILCA constructs each element of $\Im(\mathcal{R}_{\A})$ only once, reuses it whenever an identical transition submatrix recurs, and recovers scalar-equivalent RRPEFs by direct rescaling. 

\section{ILCA in Linear Programming}\label{sec. ILCA in Linear Programming} Throughout this section, $\Field$ is specialized to the real field $\reals$. Let $\C\in\reals^{\rho\times\nu}$ be a full-row-rank constraint matrix, with $\rank(\C)=\rho\leq\nu$, and let $\RS(\C)\subseteq\reals^{\nu}$ denote its row space. Consider the objective function 
\[ z=\sum_{j=1}^{\nu}z_jy_j=\z\y, \] 
where $\z=[z_j]\in\reals^{1\times\nu}$ and $\y=[y_j]\in\reals^{\nu\times1}$. The linear program is written as 
\[ \max/\min\;z=\z\y\qquad\text{s.t.}\qquad\C\y\leq\mathbf{b},\quad\y\geq\0, \] 
where $\mathbf{b}\in\reals^{\rho\times1}$ is the right-hand side. Introducing slack variables $\mathbf{s}\in\reals^{\rho\times1}$ gives \[ \x\triangleq\begin{bmatrix}\y\\\mathbf{s}\end{bmatrix},\qquad\widetilde{\z}\triangleq[\z\mid\0], \] and the augmented constraint matrix \[ \A\triangleq[\C\mid\I_{\rho}]\in\reals^{\rho\times n},\qquad n\triangleq\nu+\rho. \] Because the identity block occupies the rightmost $\rho$ columns in row order, $\A$ satisfies Definition~\ref{def. LREF} and is in LRREF. The linear program becomes \[ \max/\min\;z=\widetilde{\z}\x\qquad\text{s.t.}\qquad\A\x=\mathbf{b},\quad\x\geq\0, \] with $\widetilde{\z}\x=\z\y$. Henceforth, the tilde is suppressed and $\z\in\reals^{1\times n}$ denotes the augmented objective-coefficient row. \subsection{Nonsingular Coordinate Isolation and Basic-Solution Tableaux}\label{subsec:basic_solutions} A basic solution is determined by a BCI set $\mathsf{S}\in\mathcal{S}_{\A}$. Its basic- and nonbasic-coordinate blocks are defined by \begin{equation}\label{eq:basic_solution} \x_{\mathsf{S}}\triangleq\A_{\mathsf{S}}^{-1}\mathbf{b},\qquad\x_{\mathsf{S}'}\triangleq\0. \end{equation} The associated basic-solution tableau is obtained through nonsingular coordinate isolation by multi-index Gauss--Jordan elimination. Classical simplex pivoting proceeds locally by exchanging one basic column at a time, whereas nonsingular coordinate isolation reduces the entire transition submatrix $\A_{\mathsf{S}}$ simultaneously. Both procedures, when applied to the same BCI set $\mathsf{S}$, produce the same reduced constraint matrix $\R_{\mathsf{S}}$, the same coordinate blocks defined in \eqref{eq:basic_solution}, and hence the same basic solution defined below. 
\begin{remark}[Transition Submatrices and Basis Labels]\label{rem:identity_tableaux} \rm For every $\mathsf{S}\in\mathcal{S}_{\A}$, the corresponding $\rho\times\rho$ transition submatrix $\A_{\mathsf{S}}$ is a maximal-order nonsingular submatrix of $\A$. The transition submatrix $\A_{\mathsf{S}}$ determines the basic-coordinate block $\x_{\mathsf{S}}\in\reals^{\rho}$ in \eqref{eq:basic_solution}, whereas the BCI set $\mathsf{S}$ determines its embedding in the original coordinate positions. Consequently, if \begin{subequations}\label{eq.duplicate.transition.remark} \begin{equation}\label{eq.duplicate.transition.remark.a} \A_{\mathsf{S}}=\A_{\mathsf{T}},\qquad\mathsf{S}\neq\mathsf{T}, \end{equation} then \begin{equation}\label{eq.duplicate.transition.remark.b} \R_{\mathsf{S}}=\R_{\mathsf{T}},\qquad\x_{\mathsf{S}}=\x_{\mathsf{T}}. \end{equation} \end{subequations} Nevertheless, the corresponding basic solutions $\bar{\x}_{\mathsf{S}}$ and $\bar{\x}_{\mathsf{T}}$ in $\reals^n$ need not coincide, because the common basic-coordinate block is embedded in the distinct coordinate supports $\mathsf{S}$ and $\mathsf{T}$, respectively. Within the identical RRPEFs $\R_{\mathsf{S}}$ and $\R_{\mathsf{T}}$, the distinct BCI sets must therefore be retained as parts of the identity configurations distinguishing $(\I_{\rho},\mathsf{S})$ from $(\I_{\rho},\mathsf{T})$. Thus, the unit vectors of $\I_{\rho}$ occupy the columns of $\R_{\mathsf{S}}$, respectively $\R_{\mathsf{T}}$, indexed by $\mathsf{S}$, respectively $\mathsf{T}$. \end{remark} 
Let the linear system be given in standard form: \begin{equation}\label{eq:standard_form} \A\x=\mathbf{b},\qquad\A\in\reals^{\rho\times n},\quad\x\in\reals^n,\quad\mathbf{b}\in\reals^{\rho}. \end{equation} Assume that $\A$ has full row rank: \[ \rank(\A)=\rho. \] By Theorem~\ref{theo. Rank Minor Theorem}, $\A$ contains a nonsingular submatrix \[ \A_{\mathsf{S}}\in\reals^{\rho\times\rho} \] for some BCI set $\mathsf{S}\in\mathcal{S}_{\A}$. Let \[ \A_{\mathsf{S}'}\in\reals^{\rho\times(n-\rho)} \] contain the remaining columns, indexed by the complementary set $\mathsf{S}'$. This decomposition is a \emph{nonsingular coordinate isolation}; the columns indexed by $\mathsf{S}$ and $\mathsf{S}'$ need not be consecutive in the original matrix. Let $\Per_{\mathsf{S}}$ denote the permutation matrix that maps the basic--nonbasic block order to the original coordinate order. Since \[ \Per_{\mathsf{S}}^{-1}=\Per_{\mathsf{S}}^{\top}, \] the corresponding column partition is \begin{subequations}\label{eq.coordinate.isolation} \begin{equation}\label{eq:block_partition} \A=\bmat{\A_{\mathsf{S}}&|&\A_{\mathsf{S}'}}\Per_{\mathsf{S}}^{\top}. \end{equation} The basic solution associated with $\mathsf{S}$ is denoted by $\bar{\x}_{\mathsf{S}}$ and is defined by \begin{equation}\label{eq:complete_basic_solution} \bar{\x}_{\mathsf{S}}\triangleq\Per_{\mathsf{S}}\begin{bNiceArray}{c}\x_{\mathsf{S}}\\[3pt]\cmidrule(lr){1-1}\x_{\mathsf{S}'}\end{bNiceArray}. \end{equation} Thus, $\x_{\mathsf{S}}$ and $\x_{\mathsf{S}'}$ denote the basic- and nonbasic-coordinate blocks, respectively, whereas $\bar{\x}_{\mathsf{S}}$ denotes the basic solution in the original coordinate order. Since $\A_{\mathsf{S}}$ is nonsingular, $\A_{\mathsf{S}}^{-1}$ is the transition operator that reduces the columns indexed by $\mathsf{S}$ to the identity. Applying it to \eqref{eq:block_partition} yields the RRPEF \begin{equation}\label{eq:tableau2} \begin{aligned} \R_{\mathsf{S}} &=\A_{\mathsf{S}}^{-1}\A\\ &=\A_{\mathsf{S}}^{-1}\bmat{\A_{\mathsf{S}}&|&\A_{\mathsf{S}'}}\Per_{\mathsf{S}}^{\top}\\ &=\bmat{\I_{\rho}&|&\A_{\mathsf{S}}^{-1}\A_{\mathsf{S}'}}\Per_{\mathsf{S}}^{\top}. \end{aligned} \end{equation} By \eqref{eq:basic_solution}, Eq.~\eqref{eq:complete_basic_solution} embeds the basic-coordinate block $\x_{\mathsf{S}}=\A_{\mathsf{S}}^{-1}\mathbf{b}$ and the zero nonbasic-coordinate block in their original coordinate positions. The corresponding transformed right-hand side is \begin{equation}\label{eq:rhs} \mathrm{RHS}_{\mathsf{S}}=\x_{\mathsf{S}}. \end{equation} \end{subequations} \noindent The preceding construction is summarized in the following proposition. \begin{proposition}[Invariant Properties of the Nonsingular Isolation Principle]\label{prop:invariant_properties} For every \emph{BCI} set $\mathsf{S}\in\mathcal{S}_{\A}$, the tableau determined by \eqref{eq:block_partition}--\eqref{eq:rhs} has the following properties: \begin{enumerate}[\itshape i)] \item The columns indexed by $\mathsf{S}$ form the identity matrix: \[ (\R_{\mathsf{S}})_{\mathsf{S}}=\I_{\rho}. \] \item The basic- and nonbasic-coordinate blocks are given by \eqref{eq:basic_solution}, and the basic solution in the original coordinate order is $\bar{\x}_{\mathsf{S}}$, as defined in \eqref{eq:complete_basic_solution}. \item The transformed system is equivalent to the original system: \begin{equation}\label{eq:row_equivalence} \A\x=\mathbf{b}\quad\Longleftrightarrow\quad\A_{\mathsf{S}}^{-1}\A\x=\A_{\mathsf{S}}^{-1}\mathbf{b}. \end{equation} \item The matrix $\A_{\mathsf{S}}^{-1}$ is the invertible transition operator that transforms $\A$ into $\R_{\mathsf{S}}$, as given in \eqref{eq:tableau2}. \item If $\mathsf{S},\mathsf{T}\in\mathcal{S}_{\A}$ satisfy \eqref{eq.duplicate.transition.remark.a}, then \eqref{eq.duplicate.transition.remark.b} holds. However, the corresponding basic solutions \[ \bar{\x}_{\mathsf{S}}=\Per_{\mathsf{S}}\begin{bNiceArray}{c}\x_{\mathsf{S}}\\[3pt]\cmidrule(lr){1-1}\x_{\mathsf{S}'}\end{bNiceArray}\qquad\text{and}\qquad\bar{\x}_{\mathsf{T}}=\Per_{\mathsf{T}}\begin{bNiceArray}{c}\x_{\mathsf{T}}\\[3pt]\cmidrule(lr){1-1}\x_{\mathsf{T}'}\end{bNiceArray} \] need not coincide, since $(\I_{\rho},\mathsf{S})\neq(\I_{\rho},\mathsf{T})$ whenever $\mathsf{S}\neq\mathsf{T}$. \end{enumerate} \end{proposition} For every $\mathsf{S}\in\mathcal{S}_{\A}$, let \[ \z_{\mathsf{S}}\triangleq(z_j)_{j\in\mathsf{S}} \] denote the basic-cost block, and define the transformed objective row and current objective value by \[ \bar{\z}_{\mathsf{S}}\triangleq\z-\z_{\mathsf{S}}\R_{\mathsf{S}},\qquad\bar z_{0,\mathsf{S}}\triangleq\z_{\mathsf{S}}\x_{\mathsf{S}}. \] The associated basic-solution tableau is written in the compact form \begin{equation}\label{eq:compact_tableau} \mathcal{T}_{\mathsf{S}} \triangleq \left(\R_{\mathsf{S}}\;\middle|\;\x_{\mathsf{S}}\;;\;\bar{\z}_{\mathsf{S}},\bar z_{0,\mathsf{S}}\right) = \left(\A_{\mathsf{S}}^{-1}\A\;\middle|\;\A_{\mathsf{S}}^{-1}\mathbf{b}\;;\;\bar{\z}_{\mathsf{S}},\bar z_{0,\mathsf{S}}\right). \end{equation} \begin{theorem}[Duplicate Bases and Basis-Dependent Optimality]\label{theo. duplicate bases} Let \[ \A\in\reals^{\rho\times n},\qquad\rank(\A)=\rho,\qquad\mathbf{b}\in\reals^{\rho},\qquad\z\in\reals^{1\times n}, \] and let $\mathsf{S},\mathsf{T}\in\mathcal{S}_{\A}$ be distinct \emph{BCI} sets such that
\begin{subequations}	
	 \begin{equation}\label{eq.duplicate.transition.theorem} \A_{\mathsf{S}}=\A_{\mathsf{T}}=\Qb,\qquad\mathsf{S}\neq\mathsf{T}. \end{equation} Then: \begin{enumerate}[\itshape i)] \item \begin{equation}\label{eq.duplicate.RRPEFs.basic.blocks} \R_{\mathsf{S}}=\Qb^{-1}\A=\R_{\mathsf{T}},\qquad\x_{\mathsf{S}}=\Qb^{-1}\mathbf{b}=\x_{\mathsf{T}}. \end{equation} \item The basic solutions $\bar{\x}_{\mathsf{S}}$ and $\bar{\x}_{\mathsf{T}}$ have the same feasibility classification. In general, \begin{equation}\label{eq.distinct.basic.solutions.objective.rows} \bar{\x}_{\mathsf{S}}\neq\bar{\x}_{\mathsf{T}}\qquad\text{and/or}\qquad\bar{\z}_{\mathsf{S}}\neq\bar{\z}_{\mathsf{T}}. \end{equation} \item \begin{equation}\label{eq.difference.transformed.objective.rows} \bar{\z}_{\mathsf{S}}-\bar{\z}_{\mathsf{T}}=(\z_{\mathsf{T}}-\z_{\mathsf{S}})\R_{\mathsf{S}}. \end{equation} \item \begin{equation}\label{eq.difference.objective.values} \bar z_{0,\mathsf{S}}-\bar z_{0,\mathsf{T}}=(\z_{\mathsf{S}}-\z_{\mathsf{T}})\Qb^{-1}\mathbf{b}. \end{equation} \item The corresponding basic-solution tableaux may have distinct objective values and optimality statuses. \end{enumerate}
	\end{subequations}
	  \end{theorem} \begin{proof} By \eqref{eq.duplicate.transition.theorem}, \[ \A_{\mathsf{S}}^{-1}=\Qb^{-1}=\A_{\mathsf{T}}^{-1}, \] which gives \eqref{eq.duplicate.RRPEFs.basic.blocks}. The equality of the basic-coordinate blocks, together with the zero nonbasic-coordinate blocks, gives the same feasibility classification, whereas their embeddings in the distinct coordinate supports may yield different basic solutions. Subtraction of the definitions of $\bar{\z}_{\mathsf{S}}$ and $\bar{\z}_{\mathsf{T}}$ gives \eqref{eq.difference.transformed.objective.rows}, while subtraction of the corresponding objective values gives \eqref{eq.difference.objective.values}. These identities establish the remaining statements. \end{proof} For every $\mathsf{S}\in\mathcal{S}_{\A}$, define \begin{equation}\label{eq.transformed.objective.coefficient} \bar z_{j,\mathsf{S}}\triangleq z_j-\z_{\mathsf{S}}(\R_{\mathsf{S}})_{\bullet j},\qquad j\in\I(n),\qquad\bar{\z}_{\mathsf{S}}=\bigl(\bar z_{1,\mathsf{S}},\ldots,\bar z_{n,\mathsf{S}}\bigr). \end{equation} Under the assumptions of Theorem~\ref{theo. duplicate bases}, suppose further that $\A_{\bullet p}=\A_{\bullet q}$ with $p\neq q$, and that \[ \mathsf{S}=\{s_1,\ldots,s_{k-1},p,s_{k+1},\ldots,s_{\rho}\},\qquad\mathsf{T}=\{s_1,\ldots,s_{k-1},q,s_{k+1},\ldots,s_{\rho}\}. \] For example, in a $4\times n$ matrix $\A$, let \[ p=3,\qquad q=4,\qquad\mathsf{S}=\{1,2,3,5\},\qquad\mathsf{T}=\{1,2,4,5\}. \] Then $s_i=t_i$ for $i=1,2,4$, whereas $q=t_3=4\neq3=s_3=p$. \begin{corollary}[Single-Column Duplicate Exchange]\label{cor. duplicate column exchange} Under the assumptions of Theorem~\ref{theo. duplicate bases}, suppose further that $\A_{\bullet p}=\A_{\bullet q}$ with $p\neq q$. Then: \begin{enumerate}[\itshape i)] \item \[ \bar z_{q,\mathsf{S}}=z_q-z_p. \] \item \[ \bar z_{p,\mathsf{T}}=z_p-z_q. \] \item \begin{equation}\label{eq.opposite.exchanged.transformed.objective.coefficients} \bar z_{q,\mathsf{S}}=-\bar z_{p,\mathsf{T}}. \end{equation} \item If $z_p\neq z_q$, then \[ \bar z_{q,\mathsf{S}}\,\bar z_{p,\mathsf{T}}<0. \] \end{enumerate} \end{corollary} \begin{proof} Since $\A_{\mathsf{S}}=\A_{\mathsf{T}}$, Theorem~\ref{theo. duplicate bases} gives \[ \R_{\mathsf{S}}=\R_{\mathsf{T}}=\Qb^{-1}\A. \] Since the columns $\A_{\bullet p}$ and $\A_{\bullet q}$ coincide and occupy the $k$-th position in $\A_{\mathsf{S}}$ and $\A_{\mathsf{T}}$, respectively, \[ \Qb^{-1}\A_{\bullet p}=\Qb^{-1}\A_{\bullet q}=\e_k. \] Consequently, \[ (\R_{\mathsf{S}})_{\bullet p}=(\R_{\mathsf{S}})_{\bullet q}=(\R_{\mathsf{T}})_{\bullet p}=(\R_{\mathsf{T}})_{\bullet q}=\e_k. \] Since the $k$-th entries of $\z_{\mathsf{S}}$ and $\z_{\mathsf{T}}$ are $z_p$ and $z_q$, respectively, \eqref{eq.transformed.objective.coefficient} gives \[ \bar z_{q,\mathsf{S}}=z_q-\z_{\mathsf{S}}\e_k=z_q-z_p \] and \[ \bar z_{p,\mathsf{T}}=z_p-\z_{\mathsf{T}}\e_k=z_p-z_q. \] The remaining statements follow immediately. \end{proof} \begin{remark}\label{rem.duplicate.RRPEF.tableaux} {\rm For every $\mathsf{S}\in\mathcal{S}_{\A}$, the equivalence class \begin{equation}\label{eq.RRPEF.equivalence.class} [\mathsf{S}] = \left\{\mathsf{T}\in\mathcal{S}_{\A}\ \middle|\ \mathsf{T}\sim_{\mathcal{S}_{\A}}\mathsf{S}\right\} = \left\{\mathsf{T}\in\mathcal{S}_{\A}\ \middle|\ \R_{\mathsf{T}}=\R_{\mathsf{S}}\right\} \end{equation} 
	  contains all BCI sets that generate the same RRPEF matrix $\R_{\mathsf{S}}$. Duplicate-basis detection permits reuse of the common RRPEF and basic-coordinate block, but every BCI set in \eqref{eq.RRPEF.equivalence.class} remains in the exhaustive ILCA traversal. Its corresponding basic-solution tableau and transformed objective row must be tested separately.} \end{remark} 
  Using the set-theoretic realization \begin{equation}\label{eq.Indexed.RRPEF.family} \mathcal{R}_{\A}\triangleq\left\{(\R_{\mathsf{S}},\mathsf{S})\ \middle|\ \mathsf{S}\in\mathcal{S}_{\A}\right\}, \end{equation} together with the quotient-induced mapping \begin{equation}\label{eq.Quotient.induced.RRPEF.bijection} \widetilde{\mathcal{R}}_{\A}:\mathcal{S}_{\A}^{o}\longrightarrow\Im(\mathcal{R}_{\A}),\qquad[\mathsf{S}]\longmapsto\R_{\mathsf{S}}, \end{equation} the preceding correspondences are summarized by the following commutative diagram: 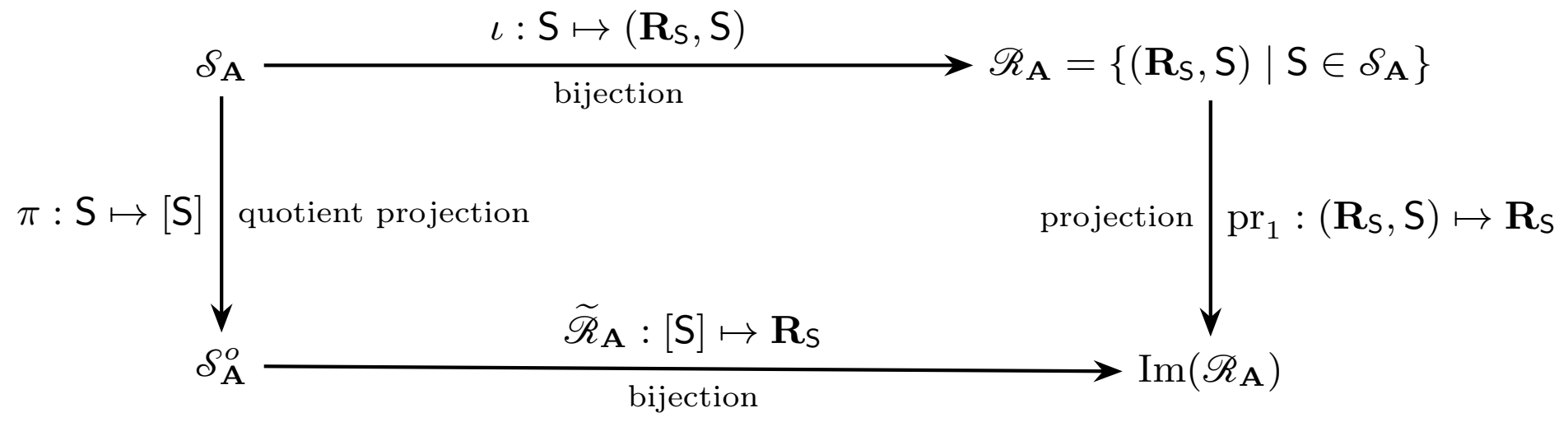
\begin{figure}[H] \begin{center} \begin{tikzpicture}[node distance=1.9cm and 5.7cm,every node/.style={font=\small,align=center},arrow/.style={-{Stealth},thick}] \node (SA) {$\mathcal{S}_{\A}$}; \node (RV) [right=of SA] {$\mathcal{R}_{\A}=\{(\R_{\mathsf{S}},\mathsf{S})\mid\mathsf{S}\in\mathcal{S}_{\A}\}$}; \node (SAo) [below=of SA] {$\mathcal{S}_{\A}^{o}$}; \node (ImRV) [below=of RV] {$\Im(\mathcal{R}_{\A})$}; \draw[arrow] (SA) -- node[above] {$\iota:\mathsf{S}\mapsto(\R_{\mathsf{S}},\mathsf{S})$} node[below] {\scriptsize bijection} (RV); \draw[arrow] (SAo) -- node[above] {$\widetilde{\mathcal{R}}_{\A}:[\mathsf{S}]\mapsto\R_{\mathsf{S}}$} node[below] {\scriptsize bijection} (ImRV); \draw[arrow] (SA) -- node[left] {$\pi:\mathsf{S}\mapsto[\mathsf{S}]$} node[right] {\scriptsize quotient projection} (SAo); \draw[arrow] (RV) -- node[right] {$\operatorname{pr}_1:(\R_{\mathsf{S}},\mathsf{S})\mapsto\R_{\mathsf{S}}$} node[left] {\scriptsize projection} (ImRV); \end{tikzpicture} \caption{Commutative correspondence between the indexed RRPEF family and its image set under quotient identification.} \label{fig:IndexedFamilyQuotientDiagram} \end{center} \end{figure} 
The upper horizontal mapping is the canonical bijection between the BCI sets and the set-theoretic realization of the indexed RRPEF family. The left vertical mapping is the quotient projection, whereas the right vertical mapping retains only the corresponding RRPEF matrix. The lower horizontal mapping is the induced bijection between the quotient set and the image set of distinct\emph{ RRPEFs}. The diagram is commutative, since, for every $\mathsf{S}\in\mathcal{S}_{\A}$, \begin{equation}\label{eq.RRPEF.diagram.commutativity} \operatorname{pr}_1\bigl(\iota(\mathsf{S})\bigr)=\R_{\mathsf{S}}=\widetilde{\mathcal{R}}_{\A}\bigl(\pi(\mathsf{S})\bigr). \end{equation} Equivalently, \[ \operatorname{pr}_1\circ\iota=\widetilde{\mathcal{R}}_{\A}\circ\pi. \] When $\R_{\mathsf{S}}=\R_{\mathsf{T}}$ with $\mathsf{S}\neq\mathsf{T}$, the indexed family retains two distinct elements and, by \eqref{eq:compact_tableau}, two distinct basic-solution tableaux may occur.

\subsection{Canonical Bijections and the Two-Level ILCA Traversal}\label{sub. bijective_identity_stalling} The ILCA traversal operates simultaneously at two synchronized levels. The RRPEF-construction level is indexed by $\mathcal{S}_{\A}^{o}$ and represented by $\Im(\mathcal{R}_{\A})$, yielding the contracted traversal \begin{equation}\label{eq.contraction.traversal} \mathcal{S}_{\A}^{o}\longleftrightarrow\Im(\mathcal{R}_{\A}). \end{equation} The basic-solution-tableau level is indexed by $\mathcal{S}_{\A}$ and contains the corresponding indexed basic-solution tableaux. Recall that the horizontal mapping \[ \mathcal{S}_{\A}\overset{\iota}{\longleftrightarrow}\mathcal{R}_{\A} \] in Figure~\ref{fig:IndexedFamilyQuotientDiagram} is the canonical bijection. The tableau construction \eqref{eq:compact_tableau} extends this correspondence. Define \begin{equation}\label{eq.Indexed.tableau.family} \mathcalligra{T}_{\A,\mathbf{b}}\triangleq\left\{(\mathcal{T}_{\mathsf{S}},\mathsf{S})\ \middle|\ \mathsf{S}\in\mathcal{S}_{\A}\right\}, \end{equation} where $\mathcal{T}_{\mathsf{S}}$ is defined by \eqref{eq:compact_tableau}. Then \begin{equation}\label{eq.RRPEF.tableau.bijection} \tau:\mathcal{R}_{\A}\longrightarrow\mathcalligra{T}_{\A,\mathbf{b}},\qquad(\R_{\mathsf{S}},\mathsf{S})\longmapsto(\mathcal{T}_{\mathsf{S}},\mathsf{S}), \end{equation} is bijective. Consequently, the exhaustive traversal is represented by \begin{equation}\label{eq.exhaustive.ILCA.traversal} \mathcal{S}_{\A}\overset{\iota}{\longleftrightarrow}\mathcal{R}_{\A}\overset{\tau}{\longleftrightarrow}\mathcalligra{T}_{\A,\mathbf{b}}. \end{equation} This level retains every BCI set and, through \eqref{eq:compact_tableau}, its basic solution and transformed objective row. At the contracted level, \eqref{eq.contraction.traversal} parametrizes the distinct RRPEF matrices, whereas \eqref{eq.exhaustive.ILCA.traversal} parametrizes the admissible BCI configurations and their indexed basic-solution tableaux. The contracted traversal constructs or recovers every distinct RRPEF in $\Im(\mathcal{R}_{\A})$, subject to the transition-contraction mechanisms. Singular-transition skipping removes candidates outside $\mathcal{S}_{\A}$; duplicate-RRPEF reuse avoids a new construction of an already obtained RRPEF; and scalar-equivalent RRPEF recovery obtains the corresponding distinct RRPEF by direct rescaling. The exhaustive traversal retains every BCI set in $\mathcal{S}_{\A}$. In the duplicate case \eqref{eq.duplicate.transition.theorem}, relation \eqref{eq.duplicate.RRPEFs.basic.blocks} holds, but, in general, the inequalities in \eqref{eq.distinct.basic.solutions.objective.rows} hold. Consequently, duplicate-RRPEF reuse removes only the repeated RRPEF construction; it removes neither BCI set nor its corresponding basic-solution tableau from the exhaustive traversal. If \begin{equation}\label{eq.scalar.equivalent.transition.submatrices} \A_{\mathsf{T}}=\lambda\A_{\mathsf{S}},\qquad\lambda\in\reals\setminus\{0,1\}, \end{equation} then \begin{equation}\label{eq.scalar.equivalent.RRPEF.basic.block} \R_{\mathsf{T}}=\lambda^{-1}\R_{\mathsf{S}},\qquad\x_{\mathsf{T}}=\lambda^{-1}\x_{\mathsf{S}}. \end{equation} Hence, the RRPEF and the basic-coordinate block are recovered by direct rescaling. If $\lambda>0$, the signs of the basic coordinates are preserved. If $\lambda<0$, the signs of all nonzero basic coordinates are reversed. In particular, if $\x_{\mathsf{S}}=\0$, then $\x_{\mathsf{T}}=\0$ and feasibility is preserved; otherwise, the feasibility classification is determined directly from the rescaled basic-coordinate block. The transformed objective row $\bar{\z}_{\mathsf{T}}$ must be tested separately. Accordingly, the transition-contraction mechanisms operate at the RRPEF-construction level, where ILCA produces each distinct RRPEF either by Gauss--Jordan reduction or through scalar-equivalent transition recovery. The resulting RRPEFs are subsequently used at the basic-solution-tableau level, yielding the exhaustive ILCA traversal, where each BCI set $\mathsf{S}\in\mathcal{S}_{\A}$ determines a distinct indexed tableau whose transformed objective row is tested separately. The strict lexicographic order on $\binom{[n]}{\rho}$ first restricts the candidate domain to \[ [\mathsf{S}_{\ell},\mathsf{S}_{\calr}]_L. \] The admissible BCI collection satisfies \[ \mathcal{S}_{\A}\subseteq[\mathsf{S}_{\ell},\mathsf{S}_{\calr}]_L, \] but equality need not hold because the interval may contain singular candidate pivot-index sets. The ILCA traversal proceeds through $\mathcal{S}_{\A}$ in strict lexicographic order, while the quotient projection $\pi$ identifies BCI sets that generate the same RRPEF. Consequently, the set-theoretic contraction is represented by \[ \binom{[n]}{\rho}\supseteq[\mathsf{S}_{\ell},\mathsf{S}_{\calr}]_L\supseteq\mathcal{S}_{\A}\xrightarrow{\pi}\mathcal{S}_{\A}^{o}. \] The first restriction implements lexicographic confinement, the second removes the singular candidates remaining in the interval, and the quotient projection identifies the BCI sets that generate the same RRPEF matrix. At the basic-solution-tableau level, every element of $\mathcal{S}_{\A}$ remains represented in $\mathcalligra{T}_{\A,\mathbf{b}}$ and undergoes its corresponding feasibility and transformed-objective-row tests. Since $\mathcal{S}_{\A}$ and $\mathcal{S}_{\A}^{o}$ are finite and strictly ordered through their lexicographic representatives, the two synchronized traversal levels terminate finitely. Upon exhaustive completion, every distinct RRPEF has been constructed or recovered, and every admissible indexed basic-solution tableau has been tested exactly once. \subsection{Structural Properties of ILCA}\label{sub. structural_integrity_convergence} The flow diagram in Figure~\ref{fig:Algorithm} is now adapted to the LP setting. The contracted traversal \eqref{eq.contraction.traversal} governs the construction of the distinct RRPEFs, whereas the exhaustive traversal \eqref{eq.exhaustive.ILCA.traversal} governs the examination of the basic-solution tableaux associated with all BCI sets in $\mathcal{S}_{\A}$. Both traversals inherit the strict lexicographic order of their corresponding BCI representatives. \begin{property}[Global Anti-Cycling under Exhaustive Traversal]\label{property anti-cycling} Let \[ \mathcal{S}_{\A}=\{\mathsf{S}_1,\ldots,\mathsf{S}_q\},\qquad\mathsf{S}_1\prec_L\cdots\prec_L\mathsf{S}_q, \] be the complete \emph{BCI} chain of $\A$. Then the exhaustive \emph{ILCA} traversal is cycle-free. Moreover, every traversal obtained from the complete \emph{BCI} chain by the transition-contraction mechanisms is cycle-free. \end{property} \begin{proof} The exhaustive traversal is strictly increasing under $\prec_L$ and therefore contains no repeated BCI set. Every contracted traversal is a subsequence of this strictly increasing chain and hence contains no repeated BCI set. Thus, neither the exhaustive nor any contracted ILCA traversal can exhibit cycling. \end{proof}

 \begin{property}[Transition-Contraction Mechanisms in the LP Setting]\label{prop.transition.contraction.LP} Let \[ \mathsf{S}_{\mathrm{next}}\in[\mathsf{S}_{\ell},\mathsf{S}_{\calr}]_L \] be the next candidate pivot-index set in the \emph{ILCA} traversal. Then: \begin{enumerate}[\itshape i)] \item \emph{Singular-transition skipping.} If \[ \det(\A_{\mathsf{S}_{\mathrm{next}}})=0, \] then \[ \mathsf{S}_{\mathrm{next}}\notin\mathcal{S}_{\A}, \] and no \emph{RRPEF} or basic-solution tableau corresponding to $\mathsf{S}_{\mathrm{next}}$ exists. \item \emph{Duplicate-RRPEF reuse.} If there exists a previously examined \emph{BCI} set $\mathsf{S}\in\mathcal{S}_{\A}$ such that \[ \A_{\mathsf{S}_{\mathrm{next}}}=\A_{\mathsf{S}}, \] then \[ \R_{\mathsf{S}_{\mathrm{next}}}=\R_{\mathsf{S}},\qquad\x_{\mathsf{S}_{\mathrm{next}}}=\x_{\mathsf{S}}. \] The common constraint matrix $\R_{\mathsf{S}_{\mathrm{next}}}=\R_{\mathsf{S}}$ is reused in the $\mathsf{S}_{\mathrm{next}}$ basic-solution tableau, whereas the transformed objective row is computed separately from the new basic-cost block: \[ \bar{\z}_{\mathsf{S}_{\mathrm{next}}}=\z-\z_{\mathsf{S}_{\mathrm{next}}}\R_{\mathsf{S}}, \] and need not coincide with $\bar{\z}_{\mathsf{S}}$. \item \emph{Scalar-equivalent transition bypassing.} If there exist a previously examined \emph{BCI} set $\mathsf{S}\in\mathcal{S}_{\A}$ and $\lambda\in\reals\setminus\{0,1\}$ such that \[ \A_{\mathsf{S}_{\mathrm{next}}}=\lambda\A_{\mathsf{S}}, \] then \[ \R_{\mathsf{S}_{\mathrm{next}}}=\lambda^{-1}\R_{\mathsf{S}},\qquad\x_{\mathsf{S}_{\mathrm{next}}}=\lambda^{-1}\x_{\mathsf{S}}. \] The basic-solution tableau corresponding to $\mathsf{S}_{\mathrm{next}}$ is retained, and its transformed objective row $\bar{\z}_{\mathsf{S}_{\mathrm{next}}}$ is tested separately. \end{enumerate} \end{property} \begin{proof} Statement~\textit{i)} follows from the BCI admissibility criterion. Statements~\textit{ii)} and~\textit{iii)} follow from Proposition~\ref{prop. transition reuse} and Eq.~\eqref{eq:basic_solution}. The separate objective-row tests follow from the exhaustive ILCA traversal \eqref{eq.exhaustive.ILCA.traversal}. \end{proof} For every BCI set $\mathsf{S}\in\mathcal{S}_{\A}$, recall that \[ \z_{\mathsf{S}}\triangleq(z_j)_{j\in\mathsf{S}} \] denotes the basic-cost block. The objective row of the transformed tableau and the current objective value are \[ \bar{\z}_{\mathsf{S}}\triangleq\z-\z_{\mathsf{S}}\R_{\mathsf{S}},\qquad\bar z_{0,\mathsf{S}}\triangleq\z_{\mathsf{S}}\x_{\mathsf{S}}. \] Accordingly, the transformed objective row augmented by the current objective value is \begin{equation}\label{eq.augmented.transformed.objective.row} \widehat{\z}_{\mathsf{S}}\triangleq\left(\bar{\z}_{\mathsf{S}}\;\middle|\;\bar z_{0,\mathsf{S}}\right). \end{equation} The implications and definitions underlying duplicate-RRPEF reuse are summarized as follows: \begin{equation}\label{eq.duplicate.RRPEF.reuse.computations} \boxed{\begin{aligned} &\mathsf{S}_{\mathrm{next}}\neq\mathsf{S}\quad\text{and}\quad\A_{\mathsf{S}_{\mathrm{next}}}=\A_{\mathsf{S}}\\ &\qquad\Longrightarrow\qquad\R_{\mathsf{S}_{\mathrm{next}}}=\R_{\mathsf{S}},\qquad\x_{\mathsf{S}_{\mathrm{next}}}=\x_{\mathsf{S}},\\[3pt] &\qquad\Longrightarrow\qquad \bar{\x}_{\mathsf{S}_{\mathrm{next}}} = \Per_{\mathsf{S}_{\mathrm{next}}} \begin{bNiceArray}{c} \x_{\mathsf{S}}\\[3pt] \cmidrule(lr){1-1} \0 \end{bNiceArray}, \qquad \bar{\x}_{\mathsf{S}} = \Per_{\mathsf{S}} \begin{bNiceArray}{c} \x_{\mathsf{S}}\\[3pt] \cmidrule(lr){1-1} \0 \end{bNiceArray},\\[5pt] &\z_{\mathsf{S}_{\mathrm{next}}}\triangleq(z_j)_{j\in\mathsf{S}_{\mathrm{next}}},\\[3pt] &\bar{\z}_{\mathsf{S}_{\mathrm{next}}}\triangleq\z-\z_{\mathsf{S}_{\mathrm{next}}}\R_{\mathsf{S}}, \qquad \bar z_{0,\mathsf{S}_{\mathrm{next}}}\triangleq\z_{\mathsf{S}_{\mathrm{next}}}\x_{\mathsf{S}},\\[3pt] &\widehat{\z}_{\mathsf{S}_{\mathrm{next}}} \triangleq \left(\bar{\z}_{\mathsf{S}_{\mathrm{next}}}\;\middle|\;\bar z_{0,\mathsf{S}_{\mathrm{next}}}\right). \end{aligned}} \end{equation} The example discussed in Appendix~\ref{append. duplicate RRPEF reuse basic-variable change} illustrates the duplicate-RRPEF reuse mechanism in Property~\ref{prop.transition.contraction.LP} and shows why the basic-solution tableau corresponding to every BCI set must be retained and tested separately. \begin{property}[Attainment of a Finite Optimum]\label{prop canonical attainment} Suppose that the \emph{LP} is feasible and has a finite optimum attained at a basic feasible solution. Then at least one basic-solution tableau examined by \emph{ILCA} encodes an optimal basic feasible solution. \end{property} \begin{proof} Every basic feasible solution is associated with at least one admissible BCI set and hence with an RRPEF basic-solution tableau. ILCA examines every admissible BCI set and therefore locates the tableau encoding an optimal basic feasible solution. Under degeneracy, more than one BCI set or tableau may encode the same optimal basic feasible solution. \end{proof} In what follows, we use the basic-solution tableau convention \[ z=\bar z_{0,\mathsf{S}}+\sum_{j\notin\mathsf{S}}\bar z_{j,\mathsf{S}}x_j\qquad\Longleftrightarrow\qquad z-\sum_{j\notin\mathsf{S}}\bar z_{j,\mathsf{S}}x_j=\bar z_{0,\mathsf{S}}. \] The objective coefficients and the current objective value are included directly in the tableau. Each Gauss--Jordan transformation therefore updates the objective row together with the constraint rows. In Tucker-type degenerate examples, the resulting tableaux may retain the same objective value throughout the RRPEF chain. \begin{property}[Terminal Halting Rule and Exhaustive Completion]\label{prop. halting_rule} \textit{i)} If $\R_{\mathsf{S}}$ is primal feasible and the objective row satisfies $\bar z_{j,\mathsf{S}}\geq0$ for all $j\notin\mathsf{S}$ in a minimization problem, respectively $\bar z_{j,\mathsf{S}}\leq0$ for all $j\notin\mathsf{S}$ in a maximization problem, then the current tableau is optimal and $z^{*}=\bar z_{0,\mathsf{S}}$. \textit{ii)} If no early halting rule is satisfied, \emph{ILCA} terminates after reaching the \emph{LRREF} and completing the finite exhaustive traversal through $\mathcal{S}_{\A}$. Provided that the \emph{LP} is feasible and has a finite optimum attained at a basic feasible solution, comparison of the objective values of the feasible basic-solution tableaux identifies the optimal value. \end{property} \begin{proof} \noindent\textit{i)} The stated feasibility and reduced-cost sign conditions are the standard sufficient optimality conditions for the adopted tableau convention. Hence, no admissible nonbasic variation can improve the current objective value. \noindent\textit{ii)} The admissible BCI collection $\mathcal{S}_{\A}$ is finite and strictly ordered. Therefore, in the absence of early termination, ILCA reaches its terminal element after constructing or recovering every distinct RRPEF and examining every basic-solution tableau associated with an admissible BCI set. By Property~\ref{prop canonical attainment}, at least one of the feasible tableaux encodes an optimal basic feasible solution. Comparison is restricted to these feasible tableaux; infeasible basic solutions do not enter the objective-value comparison. \end{proof} \noindent For linear-programming applications, the algebraic RRPEF construction is extended by adjoining the right-hand-side and objective data to each representative. The resulting basic-solution tableaux equip ILCA with the standard LP criteria for feasibility, objective evaluation, and optimality. The algorithm described below therefore admits two modes: an exhaustive traversal of the admissible BCI chain and an early-termination mode governed by the corresponding LP halting rule. 
 
 \begin{algorithmbox}[Exhaustive and Contracted ILCA Traversal] Given an LP tableau $(\A,\mathbf{b},\z)$ with $\rank(\A)=\rho$: 
 \begin{enumerate}[label=\textbf{Step \arabic*.}] 
 	\item Compute the \emph{URREF} and \emph{LRREF} of $\A$ and their corresponding \emph{BCI} sets 
 	\[ \mathsf{S}_{\ell},\qquad\mathsf{S}_{\calr}.\vspace{-0.05in} \] 
 	\item Restrict the candidate search space to 
 	\[ [\mathsf{S}_{\ell},\mathsf{S}_{\calr}]_L, \] initialize \[ \mathsf{S}\leftarrow\mathsf{S}_{\ell}, \] and form the basic-solution tableau $\mathcal{T}_{\mathsf{S}_{\ell}}$ associated with 
 	\[ \URREF(\A)=\R_{\mathsf{S}_{\ell}}. \]
 	\item Test the basic solution for feasibility and the tableau for optimality. If both conditions hold, return the corresponding optimality certificate and terminate. Otherwise, record the current objective value if the basic solution is feasible. 
 	\item If $\mathsf{S}=\mathsf{S}_{\calr}$, terminate and return the minimum or maximum objective value, as appropriate, among the feasible basic-solution tableaux examined. Otherwise, determine the next candidate 
 	\[\mathsf{S}_{\mathrm{next}}\triangleq\min_{\prec_L}\left\{\mathsf{T}\in[\mathsf{S}_{\ell},\mathsf{S}_{\calr}]_L\mid\mathsf{S}\prec_L\mathsf{T}\right\}. \vspace{-0.1in} \] 
 	\item Set 
 	\[ \mathsf{S}\leftarrow\mathsf{S}_{\mathrm{next}}. \] If $\A_{\mathsf{S}}$ is singular, return to Step~4. Otherwise, proceed to the next step. 
 	\item Form the nonsingular transition submatrix $\A_{\mathsf{S}}$. Apply duplicate-\emph{RRPEF} reuse or scalar-equivalent transition bypassing whenever applicable; otherwise, construct the corresponding \emph{RRPEF} by Gauss--Jordan reduction. In every case, form the indexed basic-solution tableau $\mathcal{T}_{\mathsf{S}}$. 
 	\item Return to Step~3. 
 \end{enumerate} \end{algorithmbox} 
\noindent Under the adopted objective-row convention, the uniform reduced-cost sign condition provides an optimality certificate for the current tableau. The strict lexicographic progression through $\mathcal{S}_{\A}$ ensures that no BCI set is revisited. Although duplicate-RRPEF reuse and scalar-equivalent transition bypassing avoid repeated Gauss--Jordan reductions, the corresponding basic-solution tableaux are nevertheless formed and tested separately. Consequently, cycling and degeneracy-induced stalling cannot obstruct finite ILCA traversal.
		
\section{ILCA on Cycling, Total Degeneracy, and Stalling Problems}\label{sec:degenerate} When applied to linear programming, ILCA augments each RRPEF with right-hand-side and objective data to form the corresponding basic-solution tableau. The resulting procedure employs the standard LP criteria for classifying feasibility, evaluating objective values, and identifying optimality. It may therefore exhaust the contracted RRPEF chain or terminate earlier when the current tableau satisfies the corresponding halting rule. All computations below are exact; the complete traces are provided in Appendices~\ref{append. beale_problem}, \ref{append.SolvingTuckerTotally DegenerateExample}, and~\ref{append.CunninghamExample}. \paragraph{Beale's cycling example.} For Beale's linear program~\cite{Beale1955}, adapted in~\cite{DantzigThapa2003}, on which Dantzig's pivot rule cycles, ILCA starts from the infeasible basic-solution tableau indexed by \[ \mathsf{S}_{\ell}=\{1,2,3\}, \] skips the three singular candidates \[ \{1,2,4\},\qquad\{1,2,5\},\qquad\{1,2,6\}, \] and proceeds through the feasible degenerate configurations $\{1,2,7\}$ and $\{1,3,4\}$. The halting rule is satisfied at \[ \mathsf{S}^{*}=\{1,3,5\},\qquad z^{*}=-\frac{1}{20}, \] which provides the optimality certificate. Thus, ILCA examines $7$ of the $35$ candidate pivot-index sets and performs $4$ RRPEF reductions; see Appendix~\ref{append. beale_problem}. \paragraph{Tucker's totally degenerate cycling example.} Consider the classical example reproduced in~\cite{DantzigThapa2003}: \[ \min\ z=-2x_1-3x_2+x_3+12x_4 \] subject to \begin{equation}\label{eq:tucker} \begin{aligned} -2x_1-9x_2+x_3+9x_4+x_5&=0,\\ \frac13x_1+x_2-\frac13x_3-2x_4+x_6&=0, \end{aligned} \qquad \x\geq\0. \end{equation} Every two-column transition submatrix is nonsingular; hence, \[ \mathcal{S}_{\A}=\binom{[6]}{2},\qquad\card(\mathcal{S}_{\A})=15. \] Moreover, the admissible BCI sets generate distinct RRPEFs. Therefore, \[ \mathcal{S}_{\A}^{o}=\mathcal{S}_{\A},\qquad\card(\mathcal{S}_{\A}^{o})=15. \] Since $\mathbf{b}=\0$, every basic solution coincides with $\x=\0$ and has objective value $z=0$. The instance is therefore totally degenerate but contains no singular or duplicate configurations. In this example, ILCA traverses the complete lexicographically ordered RRPEF chain, from the initial URREF configuration to the terminal LRREF configuration. The sequence presented in Appendix~\ref{append.SolvingTuckerTotally DegenerateExample} verifies that all $15$ configurations are generated in strict lexicographic order and that termination occurs upon exhaustion of the chain. \paragraph{Cunningham's stalling instance.} The network-simplex instance used by Cunningham to demonstrate stalling~\cite{Cunningham1979} has the initial basic-solution tableau \[ \mathcal{T}_{\{1,2,3,6,7\}}=\begin{array}{c|cccrrccrcr|c} \text{Basic}&\x_{a_0}&\x_{a_1}&\x_{a_2}&\x_{b_1}&\x_{b_2}&\x_{p_1}&\x_{p_2}&\x_{q_0}&\x_{q_1}&\x_{q_2}&\mathbf{b}\\ \hline \x_{a_0}&1&0&0&-1&0&0&0&0&0&0&0\\ \x_{a_1}&0&1&0&1&-1&0&0&0&0&0&0\\ \x_{a_2}&0&0&1&0&1&0&0&0&0&-1&0\\ \x_{p_1}&0&0&0&0&0&1&0&1&0&1&0\\ \x_{p_2}&0&0&0&0&0&0&1&-1&1&0&0\\ \hline z&0&0&0&0&0&0&0&-1&0&-1&0 \end{array}. \] Its constraint matrix is \[ \A=\R_{\{1,2,3,6,7\}}, \] with \[ \rho=5,\qquad n=10,\qquad\mathbf{b}=\0. \] Hence, all basic solutions coincide and all corresponding tableaux have the same objective value. The extreme pivot-index sets are \[ \mathsf{S}_{\ell}=\{1,2,3,6,7\},\qquad\mathsf{S}_{\calr}=\{4,5,8,9,10\}. \] The complete combinatorial layer contains \[ \card\binom{[10]}{5}=\binom{10}{5}=252 \] candidate pivot-index sets. Of these, $22$ lie outside the lexicographic interval \[ [\mathsf{S}_{\ell},\mathsf{S}_{\calr}]_L, \] with $11$ preceding $\mathsf{S}_{\ell}$ and $11$ succeeding $\mathsf{S}_{\calr}$. Among the remaining $230$ candidates, $141$ determine singular transition submatrices, leaving $89$ nonsingular BCI sets and their corresponding basic-solution tableaux. Among these $89$ nonsingular candidates, $21$ generate RRPEFs that have already been constructed because columns $7$ and $9$ coincide. For example, \[ \R_{\{1,2,3,6,9\}}=\R_{\{1,2,3,6,7\}}. \] For each such candidate, duplicate-RRPEF reuse avoids a repeated Gauss--Jordan construction of the same RRPEF, while the basic-solution tableau associated with the new BCI set is still formed and examined. Thus, ILCA examines all $89$ nonsingular basic-solution tableaux but explicitly constructs only $68$ distinct RRPEFs; see Appendix~\ref{append.CunninghamExample}. Consequently, \[ 252=22_{\mathrm{lex}}+141_{\mathrm{sing}}+21_{\mathrm{reuse}}+68_{\mathrm{constructed}}, \] and the contraction in the number of explicit RRPEF constructions is \[ 252\longrightarrow230\longrightarrow89\longrightarrow68, \] whereas the $89$ nonsingular basic-solution tableaux remain subject to feasibility and optimality testing. \paragraph{Unified view.} Cycling revisits previously encountered BCI sets or tableaux, whereas stalling proceeds through distinct tableaux without objective improvement. ILCA excludes cycling through strict lexicographic progression and remains finite independently of objective improvement. Lexicographic confinement and singular-transition skipping contract the candidate domain, whereas duplicate-RRPEF reuse reduces the number of explicit Gauss--Jordan constructions without eliminating the corresponding basic-solution tableaux. Consequently, every admissible BCI set is examined, each distinct RRPEF is constructed or recovered, and ILCA terminates either by the applicable halting rule or by exhaustion of the finite traversal.

\section*{Conclusion}\label{sec.Conclusion} This paper introduced the Reduced Row Pivot Echelon Forms (RRPEFs) as a complete canonical family of uniquely normalized, row-equivalent representatives of a matrix. Derived from a reformulation of Noether's Complementary Subspace Basis Theorem, the construction establishes a correspondence among complementary subspaces, Noether bases, admissible pivot-index sets, and RRPEFs. The URREF and LRREF form the extreme members of this family. Factoring out admissible pivot-index sets that generate identical RRPEFs yields the core combinatorial index set $\mathcal{S}_{\A}^{o}$ and the canonical bijection \[ \mathcal{S}_{\A}^{o}\longleftrightarrow\Im(\mathcal{R}_{\A}). \] This bijection transfers the strict lexicographic order to the distinct RRPEFs, forming a finite chain from the URREF to the LRREF. At the exhaustive level, the canonical correspondences \[ \mathcal{S}_{\A}\longleftrightarrow\mathcal{R}_{\A}\longleftrightarrow\mathcalligra{T}_{\A,\mathbf{b}} \] retain every admissible BCI set and its associated indexed basic-solution tableau. These two synchronized levels distinguish the construction of distinct RRPEFs from the examination of all admissible basic-solution tableaux and provide the structural link with linear programming. This ordered structure induces ILCA, which traverses the RRPEF chain in strict lexicographic order while reducing explicit RRPEF constructions through its contraction mechanisms. Feasibility and optimality are tested during traversal and may yield termination before chain exhaustion. Strict lexicographic progression excludes cycling, while the finiteness of the traversal is independent of objective improvement and therefore is not obstructed by degeneracy-induced stalling. The Beale, Tucker, and Cunningham instances illustrate complementary features of ILCA: early termination in the Beale example, complete chain traversal in the Tucker example, and substantial contraction of explicit RRPEF constructions while retaining all admissible basic-solution tableaux in the Cunningham example. Thus, the conventional single-form perspective is extended to a complete ordered family of canonical row representatives whose combinatorial structure yields an explicit finite construction procedure and an algebraic foundation for further developments of ILCA.

\section*{Future Work}
 Two principal directions emerge from the constructive framework developed in this paper. The first concerns the complexity analysis of ILCA and its extension to optimization and game-theoretic problems, including minimax and Nash-equilibrium models. A natural benchmark is the strongly polynomial result of Dadush et al.~\cite{DadushEtAl2024} for linear programs with at most two nonzero entries per row or column. The structural similarity suggest investigating whether the canonical, repetition-free RRPEF structure yields polynomial ILCA complexity on this class. 

The second direction is the \emph{Simplex--ILCA} (SILCA), developed in the companion paper~\cite{Pa:HybridLexicographicSimplex}. This hybrid method combines simplex pivoting as a local search mechanism with the finite, repetition-free RRPEF traversal of ILCA, particularly during degenerate, cycling, or stalling phases. The companion paper formulates the switching criteria, establishes global termination, and analyzes the resulting computational traversal.

\section*{Acknowledgments} The authors gratefully acknowledge Dr.~Stavros Dafnos, whose suggestion to the first author prompted the investigation of the RRPEF framework in the context of linear programming.

		\appendix
		
		\newtheorem{appexample}{Example}[section]
		\begin{center}
			\subsection*{Appendices}  
			\addcontentsline{toc}{section}{Appendices}  
		\end{center} 
		The appendices provide computational benchmarks and complements of the main body paper. 
	
\medskip	
		
\section{Matrix Singularity in Example~\ref{ex. WorkedExampleForRRPEFs}}\label{append. MatrixSingularity} This appendix provides the algebraic evaluation of an isolated coordinate singularity in Example~\ref{ex. WorkedExampleForRRPEFs} when 
\[ \det\begin{bmatrix}a&h\\b&d\end{bmatrix}=0. \] 
Setting $a=b=h=d=e=f=1$ in the matrix $\A$ of Eq.~\eqref{eq. 5x6 matrix} yields 
\[ \B=\begin{bmatrix}1&0&1&0&0&1\\0&1&1&0&0&1\\0&0&0&1&0&1\\0&0&0&0&1&1\\0&0&0&0&0&0\end{bmatrix} \] 
with 
\[ \det\begin{bmatrix}1&1\\1&1\end{bmatrix}=0. \] 
The LRREF of $\B$ is 
\[\LRREF(\B)=\begin{bmatrix}-1&1&0&0&0&0\\-1&0&-1&1&0&0\\-1&0&-1&0&1&0\\0&0&1&0&0&1\\0&0&0&0&0&0\end{bmatrix}.\vspace{0.1in} \] 
The CBI set associated with $\B$, as exhibited by $\LRREF(\B)$, is 
\[ \mathsf{C}=\{1,3\}. \] 
On the other hand, as shown by $\LRREF(\A)$ in Eq.~\eqref{equation LRREF of 5x6 matrix}, if 
\[ \det\begin{bmatrix}a&h\\b&d\end{bmatrix}\neq0, \] 
then the CBI set associated with $\A$ is $\{1,2\}$. This explains the structural discrepancy between $\LRREF(\A)$ and $\LRREF(\B)$. 

\section{Admissibility of BCI via CBI Sets}\label{append.AdmissibilityOfBCIViaCBISets} To verify admissibility through the complementary criterion, recall that a set $\mathsf{S}$ is a BCI set if and only if its complement $\mathsf{S}'$ is a CBI set. Let $\underline{\A}$ denote the matrix $\A$ of Eq.~\eqref{eq. 5x6 matrix} with its last zero row removed. A candidate CBI set $\mathsf{C}=\{i,j\}$ is admissible if the augmented determinant test gives 
\[ \det\begin{bmatrix}\underline{\A}\\\hline\e_i\\\e_j\end{bmatrix}\neq0. \] Equivalently, $\{\e_i,\e_j\}$ spans a complementary subspace of $\RS(\A)$. For the candidates $\{5,6\}$, $\{4,5\}$, and $\{4,6\}$, the test gives 
\[ \det\begin{bmatrix}\underline{\A}\\\hline\e_5\\\e_6\end{bmatrix}=0,\qquad\det\begin{bmatrix}\underline{\A}\\\hline\e_4\\\e_5\end{bmatrix}=0,\qquad\det\begin{bmatrix}\underline{\A}\\\hline\e_4\\\e_6\end{bmatrix}=0. \] Thus, none of these pairs spans a complementary subspace, and their complements are not BCI sets. Since $\A$ is in URREF, $\{3,6\}$ is an admissible CBI set, namely the complement of the BCI set $\{1,2,4,5\}$. Although the lexicographically greater candidates satisfy 
\[ \{5,6\}\succ_L\{4,6\}\succ_L\{4,5\}\succ_L\{3,6\}, \] none is admissible because each complement is lexicographically smaller than the LPP set. Applying the augmented determinant test to all candidates in $\binom{[n]}{n-\rho}$, processed in reverse lexicographic order, classifies the BCI sets through their complementary CBI sets. This provides an alternative to full-rank coordinate isolation. For $n=6$ and $\rho=4$, the $15$ pivot-index sets in $\binom{[6]}{4}$ form a lexicographic chain, and their complements form the dual chain of CBI candidates. Exactly $12$ of the $15$ complements are admissible CBI sets, forming the subchain \[ \{3,6\}\succ_L\{3,5\}\succ_L\{3,4\}\succ_L\{2,6\}\succ_L\{2,5\}\succ_L\{2,4\}\succ_L\{2,3\}\succ_L\{1,6\}\succ_L\{1,5\}\succ_L\{1,4\}\succ_L\{1,3\}\succ_L\{1,2\}. \]

\section{Duplicate-RRPEF Reuse and Separate Objective-Row Testing}\label{append. duplicate RRPEF reuse basic-variable change} Consider the linear program 
\[ \min\ z=2x_1+3x_3 \] subject to 
\[ x_1+x_2+x_3=1,\qquad x_1,x_2,x_3\geq0. \] 
Thus, 
\[ \A=\begin{bmatrix}1&1&1\end{bmatrix},\qquad\mathbf{b}=\begin{bmatrix}1\end{bmatrix},\qquad\z=\begin{bmatrix}2&0&3\end{bmatrix}. \] 
Consider the consecutive BCI sets 
\[ \mathsf{S}=\{1\},\qquad\mathsf{S}_{\mathrm{next}}=\{2\}. \] 
For $\mathsf{S}=\{1\}$, 
\[ \A_{\mathsf{S}}=\begin{bmatrix}1\end{bmatrix},\qquad\R_{\mathsf{S}}=\A_{\mathsf{S}}^{-1}\A=\begin{bmatrix}1&1&1\end{bmatrix},\qquad\x_{\mathsf{S}}=\A_{\mathsf{S}}^{-1}\mathbf{b}=\begin{bmatrix}1\end{bmatrix}. \] 
The associated basic solution is \vspace{-0.15in}
\[ \bar{\x}_{\mathsf{S}}=\begin{bmatrix}1\\0\\0\end{bmatrix}. \] 
Since 
\[ \z_{\mathsf{S}}=[z_j]_{j\in\mathsf{S}}=[z_1]=[2], \] the transformed objective row and objective value are \[ \bar{\z}_{\mathsf{S}}=\z-\z_{\mathsf{S}}\R_{\mathsf{S}}=\begin{bmatrix}2&0&3\end{bmatrix}-[2]\begin{bmatrix}1&1&1\end{bmatrix}=\begin{bmatrix}0&-2&1\end{bmatrix} \] and \[ \bar z_{0,\mathsf{S}}=\z_{\mathsf{S}}\x_{\mathsf{S}}=2. \] 
Hence, the corresponding tableau is\vspace{-0.15in} 
\[ \mathcal{T}_{\mathsf{S}}=\begin{bmatrix}x_1&x_2&x_3&\mathbf{b}\\ \hline 1&1&1&1\\\hline0&-2&1&2\end{bmatrix}. \] 
The basic solution is feasible, but the nonbasic transformed objective coefficients have opposite signs; hence, the uniform-sign optimality criterion is not satisfied. For the next BCI set, 
\[ \A_{\mathsf{S}_{\mathrm{next}}}=\begin{bmatrix}1\end{bmatrix}=\A_{\mathsf{S}}. \vspace{-0.1in} \] 
Therefore, duplicate-RRPEF reuse gives 
\[ \R_{\mathsf{S}_{\mathrm{next}}}=\R_{\mathsf{S}},\qquad\x_{\mathsf{S}_{\mathrm{next}}}=\x_{\mathsf{S}}, \] 
without a new matrix inversion or Gauss--Jordan reduction. The common basic-coordinate block is embedded in the second coordinate position, yielding 
\[ \bar{\x}_{\mathsf{S}_{\mathrm{next}}}=\begin{bmatrix}0\\1\\0\end{bmatrix}. \vspace{-0.1in}\] 
The new basic-cost block is 
\[ \z_{\mathsf{S}_{\mathrm{next}}}=[z_j]_{j\in\mathsf{S}_{\mathrm{next}}}=[z_2]=[0]. \] 
Accordingly, 
\[ \bar{\z}_{\mathsf{S}_{\mathrm{next}}}=\z-\z_{\mathsf{S}_{\mathrm{next}}}\R_{\mathsf{S}}=\begin{bmatrix}2&0&3\end{bmatrix} \] and \[ \bar z_{0,\mathsf{S}_{\mathrm{next}}}=\z_{\mathsf{S}_{\mathrm{next}}}\x_{\mathsf{S}}=0. \] The resulting tableau is 
\[ \mathcal{T}_{\mathsf{S}_{\mathrm{next}}}=\begin{bmatrix}x_1&x_2&x_3&\mathbf{b}\\\hline1&1&1&1\\\hline2&0&3&0\end{bmatrix}. \] 
Its basic solution is feasible, and its nonbasic transformed objective coefficients are nonnegative. Hence, the tableau satisfies the uniform-sign optimality criterion for minimization and yields \[ z^{*}=0. \vspace{-0.1in}\] Thus, 
\[ \R_{\mathsf{S}_{\mathrm{next}}}=\R_{\mathsf{S}},\qquad\x_{\mathsf{S}_{\mathrm{next}}}=\x_{\mathsf{S}}, \vspace{-0.1in}\] 
whereas 
\[ \bar{\x}_{\mathsf{S}_{\mathrm{next}}}\neq\bar{\x}_{\mathsf{S}},\qquad\bar{\z}_{\mathsf{S}_{\mathrm{next}}}\neq\bar{\z}_{\mathsf{S}}. \] 
Hence, duplicate-RRPEF reuse eliminates the repeated RRPEF construction but not the separate formation and optimality test of the tableau associated with $\mathsf{S}_{\mathrm{next}}$.

\section{Resolving a Beale-Type Problem via the ILCA Framework}\label{append. beale_problem} To illustrate concretely the theoretical utility of the \emph{Global Anti-Cycling Principle} in Property~\ref{property anti-cycling}, we apply ILCA to the classical degenerate configuration introduced by Beale~\cite{Beale1955} and later adapted by several authors. In the example recorded in~\cite{Yang2021}; see also Exercise~5.2 in~\cite{DantzigThapa2003}, the standard simplex method cycles indefinitely over degenerate bases: \vspace{-0.1in}
\[ \begin{aligned}\text{Minimize}\quad&z=-\frac34x_1+150x_2-\frac1{50}x_3+6x_4,\\
\text{subject to}\quad&\frac14x_1-60x_2-\frac1{25}x_3+9x_4\leq0,\\&\frac12x_1-90x_2-\frac1{50}x_3+3x_4\leq0,\\&x_3\leq1,\qquad x_1,x_2,x_3,x_4\geq0.\end{aligned} \] 
The initial canonical tableau is\vspace{-0.1in} 
\[ \mathcal{T}_{\{5,6,7\}}=\begin{bNiceArray}{ccccccc|c}[first-row,first-col,margin]&x_1&x_2&x_3&x_4&s_1&s_2&s_3&\mathbf{b}\\\hline s_1:&\frac14&-60&-\frac1{25}&9&1&0&0&0\\s_2:&\frac12&-90&-\frac1{50}&3&0&1&0&0\\s_3:&0&0&1&0&0&0&1&1\\\hline z:&-\frac34&150&-\frac1{50}&6&0&0&0&0\end{bNiceArray}. \] Applying Gauss--Jordan elimination yields the tableau \[ \mathcal{T}_{\{1,2,3\}}=\begin{bNiceArray}{ccccccc|c}[first-row,first-col,margin]&x_1&x_2&x_3&x_4&s_1&s_2&s_3&\mathbf{b}\\\hline x_1:&1&0&0&-84&-12&8&-\frac8{25}&-\frac8{25}\\x_2:&0&1&0&-\frac12&-\frac1{15}&\frac1{30}&-\frac1{500}&-\frac1{500}\\x_3:&0&0&1&0&0&0&1&1\\\hline z:&0&0&0&18&1&1&-\frac2{25}&-\frac2{25}\end{bNiceArray}, \] which is infeasible. ILCA then examines the lexicographic successors indexed by $\{1,2,4\}$, $\{1,2,5\}$, and $\{1,2,6\}$, each of which is rejected because the corresponding transition submatrix is singular. Singular-transition skipping bypasses these candidates and proceeds to the next admissible BCI set: \[ \mathcal{T}_{\{1,2,7\}}=\begin{bNiceArray}{ccccccc|c}[first-row,first-col,margin]&x_1&x_2&x_3&x_4&s_1&s_2&s_3&\mathbf{b}\\\hline x_1:&1&0&\frac8{25}&-84&-12&8&0&0\\x_2:&0&1&\frac1{500}&-\frac12&-\frac1{15}&\frac1{30}&0&0\\s_3:&0&0&1&0&0&0&1&1\\\hline z:&0&0&\frac2{25}&-18&-1&-1&0&0\end{bNiceArray}. \] The right-hand-side vector $[0,0,1]^{\top}$ is nonnegative, so the basic solution is feasible and degenerate. Since the transformed nonbasic objective coefficients are not uniformly signed, the optimality condition is not satisfied. The lexicographic traversal therefore advances to the next admissible BCI set: \[ \mathcal{T}_{\{1,3,4\}}=\begin{bNiceArray}{ccccccc|c}[first-row,first-col,margin]&x_1&x_2&x_3&x_4&s_1&s_2&s_3&\mathbf{b}\\\hline x_1:&1&-168&0&0&-\frac45&\frac{12}5&\frac2{125}&\frac2{125}\\x_3:&0&0&1&0&0&0&1&1\\x_4:&0&-2&0&1&\frac2{15}&-\frac1{15}&\frac1{250}&\frac1{250}\\\hline z:&0&-36&0&0&\frac75&-\frac{11}5&-\frac1{125}&-\frac1{125}\end{bNiceArray}. \vspace{-0.1in}\] 
Its lexicographic successor is \vspace{-0.1in}
\[ \mathcal{T}_{\{1,3,5\}}=\begin{bNiceArray}{ccccccc|c}[first-row,first-col,margin]&x_1&x_2&x_3&x_4&s_1&s_2&s_3&\mathbf{b}\\\hline x_1:&1&-180&0&6&0&2&\frac1{25}&\frac1{25}\\x_3:&0&0&1&0&0&0&1&1\\s_1:&0&-15&0&\frac{15}2&1&-\frac12&\frac3{100}&\frac3{100}\\\hline z:&0&-15&0&-\frac{21}2&0&-\frac32&-\frac1{20}&-\frac1{20}\end{bNiceArray}. \] The latter tableau is feasible because $\mathbf{b}\geq\0$ and optimal because its nonbasic transformed objective coefficients satisfy the adopted sign condition. It yields \[ z^{*}=-\frac1{20}=-0.05, \] in agreement with the classical solution. The ILCA traversal thus resolves the cyclic behavior of the standard simplex method and terminates finitely at the optimal configuration without perturbations or heuristic tie-breaking. 

\section{ILCA Applied to Tucker's Totally Degenerate Example}\label{append.SolvingTuckerTotally DegenerateExample} ILCA provides a global method for analyzing and exhaustively resolving the totally degenerate configuration associated with Tucker~\cite{Tucker1956}. Consider the simplex tableau devised by Tucker and reproduced in the book by Dantzig and Thapa~\cite{DantzigThapa2003}: 
\[ \mathcal{T}_{\{5,6\}}=\begin{bNiceArray}{cccccc|c}[first-row,first-col,margin]&x_1&x_2&x_3&x_4&s_1&s_2&\mathbf{b}\\\hline s_1:&-2&-9&1&9&1&0&0\\s_2:&\frac13&1&-\frac13&-2&0&1&0\\\hline z:&2&3&-1&-12&0&0&0\end{bNiceArray}. \] 
This tableau exhibits total degeneracy: every basic feasible solution satisfies $\mathbf{b}=\0$ and has objective value $z=0$. Under the standard simplex methodology, the pricing step selects $x_2$ as the entering variable, while the minimum-ratio test selects $s_2$ as the leaving variable. The resulting pivot is degenerate: \[ \theta=\frac{b_2}{x_{2,2}}=\frac01=0. \] Thus, the tableau changes while the objective value remains unchanged. Under the cycling pivot sequence associated with Tucker's example, successive zero-length pivots generate an infinite sequence of non-improving tableaux. Every candidate transition submatrix is nonsingular, and no two admissible BCI sets generate the same RRPEF. Therefore,\vspace{-0.1in} 
\[ \mathcal{S}_{\A}^{o}=\mathcal{S}_{\A}=\binom{[6]}{2}. \] Thus, Tucker's example is totally degenerate but not structurally singular. ILCA avoids cycling structurally. In the absence of an applicable halting rule, it follows the strictly increasing chain \begin{equation}\label{eqap. lex. chain} \R_{\{1,2\}}\prec_L\R_{\{1,3\}}\prec_L\cdots\prec_L\R_{\{4,6\}}\prec_L\R_{\{5,6\}}, \end{equation} with no possibility of reversal or repetition. Since this lexicographic chain is finite, ILCA terminates by deterministic exhaustion of $\binom{[6]}{2}$ and verifies the global value $z^{*}=0$ without cycling. In this particular example, the displayed terminal tableau can be obtained from the initial displayed tableau by a column permutation; however, column permutations do not, in general, constitute row operations and therefore do not preserve row equivalence. The initial tableau in the lexicographic ordering is \[ \mathcal{T}_{\{1,2\}}=\begin{bNiceArray}{cccccc|c}[first-row,first-col,margin]\textbf{Basic}&s_1&s_2&x_1&x_2&x_3&x_4&\mathbf{b}\\\hline s_1:&1&0&-2&-9&1&9&0\\s_2:&0&1&\frac13&1&-\frac13&-2&0\\\hline z:&0&0&2&3&-1&-12&0\end{bNiceArray}. \] The remaining tableaux in the chain~\eqref{eqap. lex. chain} are obtained by the corresponding nonsingular coordinate isolations. All $15$ tableaux are distinct valid basic representations, and their basic-variable columns form the identity matrix in their corresponding positions. Since $\mathbf{b}=\0$, every basic solution represented by these tableaux is feasible and has objective value $z=0$. Thus, despite the complete absence of singularity, every tableau is degenerate. 

\section{ILCA on Cunningham's Totally Degenerate Stalling Example}\label{append.CunninghamExample} This appendix provides the detailed combinatorial verification for Cunningham's network-simplex stalling instance discussed in the main text. Recall that $\rho=5, n=10,\mathbf{b}=\0$,  and that the extreme BCI sets are $ \mathsf{S}_{\ell}=\{1,2,3,6,7\},\mathsf{S}_{\calr}=\{4,5,8,9,10\}$.  The immediate lexicographic successor of $\mathsf{S}_{\ell}$ is $\{1,2,3,6,8\}$, whose corresponding basic-solution tableau is 
\[ \mathcal{T}_{\{1,2,3,6,8\}}=\begin{bNiceArray}{cccccccccr|c}[first-row,first-col,margin]&\x_{a_0}&\x_{a_1}&\x_{a_2}&\x_{b_1}&\x_{b_2}&\x_{p_1}&\x_{p_2}&\x_{q_0}&\x_{q_1}&\x_{q_2}&\mathbf{b}\\\hline\x_{a_0}:&1&0&0&-1&0&0&0&0&0&0&0\\\x_{a_1}:&0&1&0&1&-1&0&0&0&0&0&0\\\x_{a_2}:&0&0&1&0&1&0&0&0&0&-1&0\\\x_{p_1}:&0&0&0&0&0&1&1&0&1&1&0\\\x_{q_0}:&0&0&0&0&0&0&-1&1&-1&0&0\\\hline z:&0&0&0&0&0&0&-1&0&-1&-1&0\end{bNiceArray}. \] The $11$ candidates strictly preceding $\mathsf{S}_{\ell}$ are \[ \begin{aligned}\{1,2,3,4,5\}&\prec_L\{1,2,3,4,6\}\prec_L\{1,2,3,4,7\}\prec_L\{1,2,3,4,8\}\prec_L\{1,2,3,4,9\}\\&\prec_L\{1,2,3,4,10\}\prec_L\{1,2,3,5,6\}\prec_L\{1,2,3,5,7\}\prec_L\{1,2,3,5,8\}\\&\prec_L\{1,2,3,5,9\}\prec_L\{1,2,3,5,10\}\prec_L\mathsf{S}_{\ell}.\end{aligned} \] The $11$ candidates strictly succeeding $\mathsf{S}_{\calr}$ are \[ \begin{aligned}\mathsf{S}_{\calr}&\prec_L\{4,6,7,8,9\}\prec_L\{4,6,7,8,10\}\prec_L\{4,6,7,9,10\}\\&\prec_L\{4,6,8,9,10\}\prec_L\{4,7,8,9,10\}\prec_L\{5,6,7,8,9\}\\&\prec_L\{5,6,7,8,10\}\prec_L\{5,6,7,9,10\}\prec_L\{5,6,8,9,10\}\\&\prec_L\{5,7,8,9,10\}\prec_L\{6,7,8,9,10\}.\end{aligned} \] Thus, lexicographic confinement restricts the traversal to \[ \left|[\mathsf{S}_{\ell},\mathsf{S}_{\calr}]_L\right|=252-22=230 \] candidates. Among them, $141$ determine singular transition submatrices, leaving \[ 230-141=89 \] nonsingular BCI sets. Since columns $7$ and $9$ of $\A$ coincide, \[ \A_{\bullet7}=\A_{\bullet9}, \] replacing index $7$ by index $9$ in any listed BCI set leaves the corresponding transition submatrix unchanged, because the exchanged column occupies the same position in the increasing representation. This yields the following $21$ pairs of equal RRPEFs: \[ \resizebox{\textwidth}{!}{$\begin{alignedat}{6} \R_{\{1,2,3,6,7\}}&=\R_{\{1,2,3,6,9\}},&\qquad\R_{\{1,2,3,7,10\}}&=\R_{\{1,2,3,9,10\}},&\qquad\R_{\{1,2,5,6,7\}}&=\R_{\{1,2,5,6,9\}},\\ \R_{\{1,2,5,7,10\}}&=\R_{\{1,2,5,9,10\}},&\R_{\{1,2,6,7,10\}}&=\R_{\{1,2,6,9,10\}},&\R_{\{1,3,4,6,7\}}&=\R_{\{1,3,4,6,9\}},\\ \R_{\{1,3,4,7,10\}}&=\R_{\{1,3,4,9,10\}},&\R_{\{1,3,5,6,7\}}&=\R_{\{1,3,5,6,9\}},&\R_{\{1,3,5,7,10\}}&=\R_{\{1,3,5,9,10\}},\\ \R_{\{1,4,5,6,7\}}&=\R_{\{1,4,5,6,9\}},&\R_{\{1,4,5,7,10\}}&=\R_{\{1,4,5,9,10\}},&\R_{\{1,4,6,7,10\}}&=\R_{\{1,4,6,9,10\}},\\ \R_{\{1,5,6,7,10\}}&=\R_{\{1,5,6,9,10\}},&\R_{\{2,3,4,6,7\}}&=\R_{\{2,3,4,6,9\}},&\R_{\{2,3,4,7,10\}}&=\R_{\{2,3,4,9,10\}},\\ \R_{\{2,4,5,6,7\}}&=\R_{\{2,4,5,6,9\}},&\R_{\{2,4,5,7,10\}}&=\R_{\{2,4,5,9,10\}},&\R_{\{2,4,6,7,10\}}&=\R_{\{2,4,6,9,10\}},\\ \R_{\{3,4,5,6,7\}}&=\R_{\{3,4,5,6,9\}},&\R_{\{3,4,5,7,10\}}&=\R_{\{3,4,5,9,10\}},&\R_{\{4,5,6,7,10\}}&=\R_{\{4,5,6,9,10\}}. \end{alignedat}$} \] These identities account for the $21$ duplicate-RRPEF reuse cases among the $89$ nonsingular BCI sets. For each pair, ILCA avoids repeating the Gauss--Jordan construction of the common RRPEF while still forming and examining the basic-solution tableau associated with each BCI set. Hence, all $89$ nonsingular basic-solution tableaux are examined, whereas only \[ 89-21=68 \] distinct RRPEFs require explicit Gauss--Jordan construction. The exhaustive traversal terminates at the basic-solution tableau $\mathcal{T}_{\{4,5,8,9,10\}}$, whose constraint matrix is the LRREF \[ \R_{\{4,5,8,9,10\}}=\LRREF(\A). \] Since $\mathbf{b}=\0$, every basic solution is the zero vector and every corresponding tableau has objective value zero. Thus, the exhaustive traversal examines all $89$ nonsingular BCI sets and their basic-solution tableaux, while duplicate-RRPEF reuse reduces the number of explicit Gauss--Jordan constructions to $68$.

\end{document}